\documentstyle[12pt,a4]{article}
\title{The Mumford Relations and the Moduli of Rank Three Stable Bundles}
\author{Richard Earl,\\Mathematical Institute,\\24-29 St.
Giles,\\Oxford, OX1 3LB,\\England}
\date{}
\begin{document}
\setcounter{page}{1}
\newtheorem{REM}{Remark}
\newenvironment{rem}{\begin{REM} \rm}{\end{REM}}
\newtheorem{prop}{PROPOSITION}
\newtheorem{lem}[prop]{LEMMA}
\newtheorem{cor}[prop]{COROLLARY}
\newtheorem{thm}[prop]{THEOREM}
\newtheorem{guess}{CONJECTURE}
\maketitle
\newcommand{\G}{{\cal G}}
\newcommand{\pa}{\overline{\partial}}
\newcommand{\g}{\bar{g}}
\newcommand{\Cmu}{{\cal C}_{\mu}}
\newcommand{\om}{\tilde{\Omega}}
\newcommand{\HG}{H^{*}_{\cal G}}
\newcommand{\ch}{\mbox{ch}}
\newcommand{\dai}{\frac{\partial}{\partial a_{i}}}
\newcommand{\daij}{\frac{\partial^{2}}{\partial a_{i} \partial a_{j}}}
\newcommand{\daii}{\frac{\partial}{\partial a_{i}^{I}}}
\newcommand{\dajj}{\frac{\partial}{\partial a_{j}^{J}}}
\newcommand{\Css}{{\cal C}^{ss}}
\newcommand{\Q}{{\bf Q}}
\newcommand{\Mnd}{{\cal M}(n,d)}
\newcommand{\mnd}{{\cal M}(n,d)}
\newcommand{\HS}{H^{*}}
\newcommand{\ar}{a_{r}}
\newcommand{\bjr}{b_{r}^{s}}
\newcommand{\fr}{f_{r}}
\newcommand{\C}{{\cal C}}
\newcommand{\Cnu}{{\cal C}_{\nu}}
\newcommand{\Z}{{\bf Z}}
\newcommand{\ta}{\tilde{a}}
\newcommand{\tb}{\tilde{b}}
\newcommand{\tf}{\tilde{f}}
\newcommand{\tta}{\hat{a}}
\newcommand{\ttb}{\hat{b}}
\newcommand{\ttf}{\hat{f}}
\begin{abstract}
\indent The cohomology ring of the moduli space $\Mnd$ of semistable
bundles of coprime rank $n$ and degree $d$ over a Riemann surface $M$
of genus $g \geq 2$ has again proven a rich source of interest in
recent years. The rank two, odd degree case is now largely understood.
In 1991 Kirwan \cite{K2} proved two long standing conjectures due to Mumford
and to Newstead and Ramanan. Mumford conjectured that a certain set of relations form
a complete set; the Newstead-Ramanan conjecture involved the vanishing of the
Pontryagin ring. The Newstead-Ramanan conjecture was independently proven by
Thaddeus \cite{T} as a corollary to determining the intersection pairings.\\
\indent As yet though, little work has been done on the cohomology ring
in higher rank cases. A simple numerical calculation shows that the
Mumford relations themselves are not generally complete when $n>2$.
However by generalising the methods of \cite{K2} and by
introducing new relations, in a sense dual to the original relations
conjectured by Mumford, we prove results corresponding to the Mumford
and Newstead-Ramanan conjectures in the rank three case. Namely we show
($\S4$) that the Mumford relations and these `dual' Mumford
relations form a complete set for the rational cohomology ring of
${\cal M}(3,d)$ and show ($\S5$) that the Pontryagin ring vanishes in degree
$12g-8$ and above.
\end{abstract}
{\bf Acknowledgements:} I would like to thank Dr. Frances
Kirwan and Luke Hatter for their time and comments on this article. I
would further like to thank Dr. Alastair King and Dr. Peter Newstead
for informing me of their and other recent results. 
\section{Introduction.}
Let $\mnd$ denote the moduli space of semistable holomorphic
vector bundles of coprime rank $n$ and degree $d$ over a Riemann
surface $M$ of genus $g \geq 2.$ Throughout this article we will write
\[
\g = g-1.
\]
Recall that a holomorphic vector
bundle $E$ over $M$ is said to be semistable (resp. stable) if every proper
subbundle $F$ of $E$ satisfies
\[
\mu(F) \leq \mu(E) \indent (\mbox{resp. } \mu(F) < \mu(E))
\]
where $\mu(F) = \mbox{degree}(F)/\mbox{rank}(F)$ is the slope of $F$.
Non-semistable bundles are said to be unstable. When $n$ and $d$ are
coprime the stable and semistable bundles coincide.\\
\indent Let ${\cal E}$ be a fixed $C^{\infty}$ complex vector bundle
of rank $n$ and degree $d$ over $M$. Let $\C$ be the space of all holomorphic
structures on ${\cal E}$ and let $\G_{c}$ denote the group of all
$C^{\infty}$ complex automorphisms of ${\cal E}$. Atiyah and Bott
\cite{AB} identify the moduli space
$\mnd$ with the quotient $\Css/{\cal G}_{c}$ where $\Css$ is the open
subset of $\C$ consisting of all semistable holomorphic structures on
${\cal E}$. In this construction both $\C$
and $\G_{c}$ are infinite dimensional; there
exist other constructions \cite{K3} of the moduli space $\mnd$ as
genuine geometric invariant theoretic quotients which are in a sense
finite dimensional approximations of Atiyah and Bott's construction.\\
\indent There is a  known set of generators \cite{N2,AB} for the
rational cohomology ring of $\mnd$ as follows. Let $V$ denote a
universal bundle over $\mnd \times M$. Atiyah and Bott then define
elements
\begin{equation}
\ar \in H^{2r}(\mnd;\Q), \quad \bjr \in
H^{2r-1}(\mnd;\Q), \quad \fr \in H^{2r-2}(\mnd;\Q) \label{000}
\end{equation}
where $1 \leq r \leq n,1 \leq s \leq 2g$ by writing
\begin{equation}
c_{r}(V) = \ar \otimes 1 + \sum_{s=1}^{2g} \bjr \otimes
\alpha_{s} + \fr \otimes \omega \indent 1 \leq r \leq n \label{0}
\end{equation}
where $\omega$ is the standard generator of $H^{2}(M;\Q)$ and
$\alpha_{1},...,\alpha_{2g}$ form a fixed canonical cohomology basis for
$H^{1}(M;\Q)$. The ring
$H^{*}(\mnd;\Q)$ is freely generated as a graded algebra over $\Q$ by
the elements (\ref{000}). Notice from the definition that $f_{1}=d$. We further introduce the notation
\[
\xi_{i,j} = \sum_{s=1}^{g} b_{i}^{s} b_{j}^{s+g}.
\]
\indent The universal bundle $V$ is not unique, although its
projective class is. We may tensor $V$ by the pullback to $\mnd \times
M$ of any holomorphic line bundle $K$ over $\Mnd$ to give another
bundle with the same universal property. This process changes the
generators of $\HS(\Mnd;\Q)$. In particular it changes $a_{1}$ by $n
c_{1}(K)$ and $c_{1}(\pi_{!}V)$ by $(d-n\g)c_{1}(K)$ where $\pi:\mnd
\times M \rightarrow \mnd$ is the first projection and $\pi_{!}$ is
the direct image map from K-theory \cite[p.436]{H}. Since $n$ and $d$
are coprime there exist integers $u$ and $v$ such that
\[
u n + v (d-n\g) =1.
\]
Thus if we take $K$ to be
\[
\left. \mbox{det}(V \right|_{\Mnd}) ^{u} \otimes (\mbox{det}
\pi_{!}V)^{v}
\]
then $V \otimes \pi^{*}(K^{-1})$ is a new universal bundle such that
\begin{equation}
 u a_{1} + v c_{1}(\pi_{!}V) = 0. \label{NORM}
\end{equation}
Following Atiyah and Bott \cite[p.582]{AB} we replace $V$ by this
normalised universal bundle.\\
\indent  The normalised bundle $V$ is universal in
the sense that its
restriction to $\{[E]\} \times M$ is isomorphic to $E$ for each
semistable holomorphic bundle $E$ over $M$ of rank $n$ and degree $d$
and where $[E]$ is the class of $E$ in $\Mnd.$ Then the stalk of the $i$th higher
direct image sheaf $R^{i}\pi_{*}V$ (see \cite[$\S 3.8$]{H}) at $[E]$
is
\[
H^{i}(\pi^{-1}([E]),V_{|\pi^{-1}([E])}) = H^{i}(M,V_{|[E] \times M}) \cong
H^{i}(M,E).
\]
\indent Tensoring $E$ with a holomorphic line bundle over $M$ of
degree $D$ gives an isomorphism between $\Mnd$ and ${\cal M}(n,d+nD)$.
Since $n$ and $d$ are coprime we may assume without any loss of
generality that $2\g n<d<(2\g+1)n$ and so we will write
\[
d=2n\g+\delta \indent (0 < \delta < n)
\]
from now on. From \cite[lemma 5.2]{N} we know that $H^{1}(M,E)=0$
for any semistable holomorphic bundle $E$ of slope greater than
$2\g$. Thus $\pi_{!}V$ is in fact a vector bundle over $\Mnd$ with fibre
$H^{0}(M,E)$ over $[E] \in \Mnd$ and, by the Riemann-Roch theorem, of
rank $d-n\g=n\g+\delta$.\\
\indent In particular if we express the Chern classes
$c_{r}(\pi_{!}V)$ in terms of the generators $a_{r},b_{r}^{s}$ and
$f_{r}$ of $\HS(\mnd;\Q)$ then knowing the images of the $r$th Chern
classes in
$H^{*}(\Mnd;\Q)$ vanish for $r>n\g+\delta$ gives us relations in terms of
the images of the generators  in $H^{*}(\Mnd;\Q).$ Now from
\cite[prop. 9.7]{AB} we know that
\begin{equation}
H^{*}(\Mnd;\Q) \cong H^{*}({\cal M}_{0}(n,d);\Q) \otimes
H^{*}(\mbox{Jac}(M);\Q) \label{21}
\end{equation}
where $\mbox{Jac}(M)$ is the Jacobian of the Riemann surface $M$ and
${\cal M}_{0}(n,d)$ is the moduli space of rank $n$ bundles with
degree $d$ and fixed determinant line bundle.
$H^{*}(\mbox{Jac}(M);\Q)$ is an exterior algebra on $2g$ generators
and we can choose the isomorphism (\ref{21}) so that these generators
correspond to $b_{1}^{1},...,b_{1}^{2g}$ and the elements
$a_{2},...,a_{n},b_{2}^{1},...,b_{n}^{2g},f_{2},...,f_{n}$ correspond
to the generators of $H^{*}({\cal M}_{0}(n,d);\Q)$. So we can find
relations in terms of
$a_{2},...,a_{n},b_{2}^{1},...,b_{n}^{2g},$ and $f_{2},...,f_{n}$ by equating
to zero the coefficients of $\prod_{s \in S}b_{1}^{s}$ in the Chern
classes $c_{r}(\pi_{!}V)$ for $r>n\g+\delta$ and for every subset $S
\subseteq \{1,...,2g\}.$\\
\indent  Mumford's conjecture, as proven by Kirwan \cite[$\S$2]{K2}, was that when the rank $n$ is two then these
relations
together with the relation (\ref{NORM}) from
normalising the universal bundle $V$ provide a complete set of relations in
$H^{*}({\cal M}_{0}(2,d);\Q)$. Subsequently a stronger version of
Mumford's conjecture has been proven \cite{E} showing the relations
coming from the first vanishing Chern class $c_{2g}(\pi_{!}V)$ generate the
relation ideal of $\HS({\cal M}_{0}(2,d))$ as a $\Q[a_{2},f_{2}]$-module.\\
\begin{rem}
In the rank two case the Mumford relations above
differ somewhat from the relations $\xi_{r}$ introduced by Zagier and
studied in \cite{B,KN,ST,Z}. In the notation of \cite{Z}
\[
\Psi_{\{1,...,2g\}} \left( \frac{-t-a_{1}}{2} \right) =
\frac{(-1)^{g\g/2 + g}}{2^{2g-1}} t^{\g} F_{0}(t^{-1})
\]
where $\Psi_{\{1,...2g\}}(x)$ denotes the coefficient of
$\prod_{s=1}^{2g} b_{1}^{s}$ in $\Psi(x) =\sum_{r \geq 0}
c_{r}(\pi_{!}V) x^{2g-1-r}$ and $F_{0}(t) =
\sum_{r=0}^{\infty} \xi_{r} t^{r}$. In the notation of \cite{KN}
$\xi_{r}$ appears as $\zeta_{r}/r!$ and in \cite{ST} as
$\Phi^{(r)}/r!.$
\end{rem}
\indent We will demonstrate later (remark \ref{inadequacy}) that the Mumford relations are
not complete when the rank $n$ is greater than two. For now we
introduce a new set of relations. Let $L$ be a fixed line bundle over
$M$ of degree $4\g+1$ and 
let $\phi:\mnd \times M \rightarrow M$ be the second
projection. Then $\pi_{!}(V^{*} \otimes
\phi^{*}L)$ is a vector bundle over $\Mnd$ of rank $(3\g+1)n-d
=ng-\delta$ with fibre $H^{0}(M,E^{*} \otimes L)$ over $[E]$. By
equating to zero the coefficients of $\prod_{s \in S}b_{1}^{s}$ in the Chern
classes $c_{r}(\pi_{!}(V^{*} \otimes \phi^{*}L))$ for $r>ng- \delta$
and for every subset $S \subseteq \{1,...,2g\}$ we may find
relations in terms of the generators
$a_{2},...,a_{n},b_{2}^{1},...,b_{n}^{2g},$ and $f_{2},...,f_{n}$. We
will refer to these new relations as the dual Mumford relations.
\begin{rem}
\label{dualise}
The map $E
\mapsto E^{*} \otimes L$ induces an automorphism of $\HS({\cal
M}(2,d);\Q)$ mapping the Mumford relations to the dual Mumford
relations and vice versa. Hence we can deduce that the dual Mumford relations
are complete when the rank is two from Kirwan's proof of Mumford's
conjecture \cite[$\S$ 2]{K2}.
\end{rem}
\indent Our first result (to be proved in $\S 4$) now reads
as:\\[\baselineskip]
{\bf THEOREM 1.} {\em The Mumford and dual Mumford relations together with the
relation (\ref{NORM}) due to the normalisation of the
universal bundle $V$ form a complete set of relations for $H^{*}({\cal
M}(3,d);\Q).$}\\[\baselineskip]
\indent The Newstead-Ramanan conjecture states \cite[$\S$5a]{N2} that the
Pontryagin ring of the tangent bundle to ${\cal M}(2,d)$
vanishes in degrees $4g$ and higher. The conjecture was proven
independently by Thaddeus \cite{T} and Kirwan \cite[$\S$4]{K2}, and
has been proven more recently by King and Newstead \cite{KN} and Weitsman \cite{W}. In $\S
5$ we will use a similar method to Kirwan's but now also involving the
dual Mumford relations to prove:\\[\baselineskip]
{\bf THEOREM 2.} {\em The Pontryagin ring of the moduli space ${\cal
M}(3,d)$ vanishes in degrees $12g-8$ and above.}
\section{Kirwan's Approach.}
\indent The group $\G_{c}$ is the
complexification of the gauge group $\G$ of all smooth automorphisms
of ${\cal E}$ which are unitary with respect to a fixed Hermitian
structure on ${\cal E}$ \cite[p.570]{AB}. We shall write
$\overline{\G}$ for the quotient of $\G$ by its $U(1)$-centre and
$\overline{\G}_{c}$ for the quotient of $\G_{c}$ by its ${\bf
C}^{*}$-centre. \\
\indent There are natural isomorphisms \cite[9.1]{AB}
\[
H^{*}(\Css/\G_{c};\Q) = H^{*}(\Css/\overline{\G}_{c};\Q) \cong
H^{*}_{\overline{\G}_{c}}(\Css;\Q) \cong H^{*}_{\overline{\G}}(\Css;\Q)
\]
since the ${\bf C}^{*}$-centre of $\G_{c}$ acts trivially on $\Css$,
$\overline{\G}_{c}$ acts freely on $\Css$ and $\overline{\G}_{c}$ is
the complexification of $\overline{\G}$. Atiyah and Bott \cite[thm.
7.14]{AB} show that the restriction map $H^{*}_{\overline{\G}}(\C;\Q)
\rightarrow H^{*}_{\overline{\G}}(\Css;\Q)$ is surjective. Further
$H^{*}_{\overline{\G}}(\C;\Q) \cong
H^{*}(B\overline{\G};\Q)$ since $\C$ is an affine space \cite[p.565]{AB}. So
putting this all together we have
\begin{equation}
H^{*}(B\overline{\G};\Q) \cong H^{*}_{\overline{\G}}(\C;\Q)
\rightarrow H^{*}_{\overline{\G}}(\Css;\Q) \cong H^{*}(\mnd;\Q) \label{14}
\end{equation}
is a surjection.\\
\indent As shown in \cite[prop. 2.4]{AB} the classifying space $B\G$
can be identified with the space $\mbox{Map}_{d}(M,BU(n))$ of all
smooth maps $f:M \rightarrow BU(n)$ such that the pullback to $M$ of
the universal vector bundle over $BU(n)$ has degree $d$. If we
pull back this universal bundle using the evaluation map
\[
\mbox{Map}_{d}(M,BU(n)) \times M \rightarrow BU(n): (f,m) \mapsto f(m)
\]
then we obtain a rank $n$ vector bundle ${\cal V}$ over $B\G \times
M$. If we restrict the pullback bundle induced by the maps
\[
\Css \times E\G \times M \rightarrow \C \times E\G \times M \rightarrow \C
\times_{\G} E\G \times M \stackrel{\simeq}{\rightarrow} B\G \times M
\]
to $\Css \times \{e\} \times M$ for some $e \in E\G$ then we obtain
a $\G$-equivariant holomorphic bundle on $\Css \times M$. The $U(1)$-centre of
$\G$ acts as scalar multiplication on the fibres,
and the associated projective bundle descends to a holomorphic
projective bundle over $\mnd \times M$ which is in fact the projective
bundle of $V$ \cite[pp.579-580]{AB}.\\
\indent By a slight abuse of notation we define elements $\ar, \bjr, \fr$
in $\HS(B\G;\Q)$ by writing
\[
c_{r}({\cal V}) = \ar \otimes 1 + \sum_{s=1}^{2g} \bjr \otimes
\alpha_{s} + \fr \otimes \omega \indent 1 \leq r \leq n.
\]
Atiyah and Bott show \cite[prop. 2.20]{AB} that the ring
$H^{*}(B\G;\Q)$ is freely generated 
as a graded algebra over $\Q$ by the elements $\ar, \bjr, \fr$. The
only relations amongst these generators are that the 
$\ar$ and $\fr$ commute with everything else and that the
$\bjr$ anticommute with each other.\\
\indent The fibration $BU(1) \rightarrow B\G \rightarrow
B\overline{\cal G}$ induces an isomorphism \cite[p.577]{AB}
\[
H^{*}(B\G;\Q) \cong H^{*}(B\overline{\G};\Q) \otimes H^{*}(BU(1);\Q).
\]
The generators $\ar,\bjr$ and $\fr$ of $\HS(B\G;\Q)$ can be pulled
back via a section of this 
fibration to give rational generators of the cohomology ring of
$B\overline{\cal G}$. We may if we wish omit $a_{1}$ since its image
in $H^{*}(B\overline{\G};\Q)$ can be expressed in terms of the other
generators. The only other relations are again the commuting of the $\ar$ and
$\fr$, and the anticommuting of the $\bjr$. We may then normalise ${\cal V}$
suitably so that these generators for $H^{*}(B\overline{\G};\Q)$
restrict to the generators $\ar, \bjr, \fr$ for $\HS(\mnd;\Q)$ under
the surjection (\ref{14}).\\
\indent The relations amongst these generators for $\HS(\mnd;\Q)$ are
then given by the kernel of the restriction map (\ref{14}) which in
turn is determined by the map
\[
\HG(\C;\Q) \cong  H^{*}_{\overline{\G}}(\C;\Q) \otimes
\HS(BU(1);\Q) \rightarrow \indent\indent\indent\indent\indent\indent
\]
\[
\indent\indent\indent\indent\indent
H^{*}_{\overline{\G}}(\Css;\Q) \otimes
\HS(BU(1);\Q) \cong \HG(\Css;\Q).
\]
In order to describe this kernel we consider Shatz's stratification of
$\C$, the space of holomorphic structures on ${\cal E}$
\cite{Sh}. The stratification $\{\Cmu : \mu \in {\cal M} \}$ is
indexed by the partially ordered set ${\cal M}$, consisting of all the
types of holomorphic bundles of rank $n$ and degree $d$, as follows.\\
\indent Any holomorphic bundle $E$ over $M$ of rank $n$ and degree $d$ has a
canonical filtration (or flag) \cite[p.221]{HN}
\[
0 = E_{0} \subset E_{1} \subset \cdot \cdot \cdot \subset E_{P} = E
\]
of sub-bundles such that the quotient bundles $Q_{p}= E_{p}/E_{p-1}$
are semi-stable and $\mu(Q_{p}) > \mu(Q_{p+1})$. We
will write $d_{p}$ and $n_{p}$ respectively for the degree and rank of
$Q_{p}$. Given such a filtration we define the type of $E$ to be
\[
\mu = (\mu(Q_{1}),...,\mu(Q_{P})) \in \Q^{n}
\]
where the entry $\mu(Q_{p})$ is repeated $n_{p}$ times. When
there is no chance of confusion we will
also refer collectively to the strata of type $(n_{1},...,n_{s})$ and
we will write $\Delta$ for the collection of strata with $n_{p}=1$ for each
$p$. The
semistable bundles have type $\mu_{0} = (d/n,...,d/n)$ and form the
unique open stratum. The set ${\cal M}$ of all possible types of holomorphic
vector bundles over $M$ will provide our indexing set. A partial order on
${\cal M}$ is defined as follows. Let $\sigma=(\sigma_{1},...,\sigma_{n})$ and
$\tau=(\tau_{1},...,\tau_
{n})$ be two types; we say that $\sigma \geq \tau$ if and only if
\[
\sum_{j \leq i} \sigma_{j} \geq \sum_{j \leq i} \tau_{j} \mbox{ for } 1 \leq i
\leq n-1.
\]
The set $\Cmu \subseteq {\cal C},$ $\mu \in {\cal M}$, is
defined to be the set of all holomorphic vector bundles of type $\mu$.\\
\indent The stratification also has the following properties:-\\
\indent (i) The stratification is smooth. That is each stratum $\Cmu$ is a
locally closed $\G_{c}$-invariant
submanifold. Further for any $\mu \in {\cal M}$ \cite[7.8]{AB}\\
\begin{equation}
\overline{\Cmu} \subseteq \bigcup_{\nu \geq \mu} \C_{\nu}. \label{11}
\end{equation}
\indent (ii) Each stratum $\Cmu$ is connected and has finite (complex)
codimension $d_{\mu}$ in $\C$. Moreover given any integer $N$ there are only
finitely many $\mu \in {\cal M}$ such that $d_{\mu} \leq N$. Further $d_{\mu}$
is given by the formula \cite[7.16]{AB}
\begin{equation}
d_{\mu}= \sum_{i>j} (n_{i}d_{j}-n_{j}d_{i}+n_{i}n_{j}\g) \label{12}
\end{equation}
where $d_{k}$ and $n_{k}$ are the degree and rank, respectively, of $Q_{k}$.\\
\indent (iii) The gauge group $\G$ acts on $\C$ preserving the
stratification which is equivariantly
perfect with respect to this action \cite[thm. 7.14]{AB}. In
particular there is an isomorphism of vector spaces
\[
H^{k}_{\G}(\C;\Q) \cong \bigoplus_{\mu \in {\cal M}}
H^{k-2d_{\mu}}_{\G}(\Cmu;\Q) = H^{k}_{\G}(\Css;\Q) \oplus \bigoplus_{\mu \neq
\mu_{0}} H^{k-2d_{\mu}}_{\G} (\Cmu;\Q).
\]
The restriction map $H^{*}_{\G}(\C;\Q) \rightarrow H^{*}_{\G}(\Css;\Q)$ is the
projection onto the summand $H^{*}_{\G}(\Css;\Q)$ and so the kernel is
isomorphic as a vector space to
\begin{equation}
\bigoplus_{k \geq 0} \bigoplus_{\mu \neq \mu_{0}} H^{k-2d_{\mu}}_{\G}
(\Cmu;\Q). \label{13}
\end{equation}
\begin{rem}
\label{inadequacy}
We can at this point use a dimension argument to
show that the Mumford relations are generally not complete when the
rank $n$ is greater than two. From the isomorphism (\ref{13}) we can see that
for the Mumford relations to be complete it is necessary that the
least degree of a Mumford relation must be less than or equal to the
smallest real codimension of an unstable stratum. The degree of
$\sigma_{r,S}^{k}$ equals $2(n\g+\delta -nr -k) -|S|$ which is least
when $r=-1,k=n-1,$ and $S=\{1,...,2g\}.$ So the smallest degree of a
Mumford relation is $2(\delta + (n-1)\g).$ However a simple
calculation minimising the codimension formula (\ref{12}) shows that the
least real codimension of an unstable stratum is $2(\delta + (n-1)\g)$
when $\delta < n/2$ and is $2(n-\delta +(n-1)\g)$ when $\delta>n/2$.
Hence the Mumford relations are not complete when $n \geq 3$ and
$\delta>n/2.$ A similar argument shows that
the dual Mumford relations are not complete when $\delta<n/2$ since
the smallest degree of a dual Mumford relation
is $2(n-\delta +(n-1)\g).$ Clearly however this simple argument does
not tell us anything concerning the union of the Mumford and dual
Mumford relations.
\end{rem}
\indent To conclude this section we will describe a set of criteria for the
completeness of a set of relations in $\HS(\mnd;\Q)$ and reformulate
the Mumford and dual Mumford relations in a way more suited to these 
criteria. Consider the formal power series
\[
c(\pi_{!}{\cal V})(t) = \sum_{r \geq 0} c_{r}(\pi_{!}{\cal V}) \cdot
t^{r} \in \HG(\C;\Q)[[t]].
\]
The vanishing of the image of $c_{r}(\pi_{!}{\cal V})$ in $H^{*}(\Mnd;\Q)$
for $r > n\g+\delta$ is equivalent to the image of $c(\pi_{!}{\cal V})(t)$ being a
polynomial of degree at most $n\g+\delta$ or equally to the image of
\[
\Psi(t) = t^{n\g+\delta}c(\pi_{!}{\cal V})(t^{-1})
\]
being a polynomial of degree at most $n\g+\delta$ in
$H^{*}(\Mnd;\Q)[t].$ If we write $\Psi(t)$ as the series
\[
\Psi(t)=\sum_{r=-\infty}^{\g} (\sigma_{r}^{0} + \sigma_{r}^{1}t +
\cdot \cdot \cdot + \sigma_{r}^{n-1}t^{n-1} ) (\om(t))^{r}
\]
where $\om(t)= t^{n}+a_{1}t^{n-1}+ \cdot \cdot \cdot +a_{n}$ then the
Mumford relations are equivalent to the vanishing of the images of
$\sigma_{r,S}^{k} (r<0,0 \leq k \leq n-1,S \subseteq \{1,...,2g\})$ in
$H^{*}({\cal M}_{0}(n,d);\Q)$ when we write
\begin{equation}
\sigma_{r}^{k} = \sum_{S \subseteq \{1,...,2g\}} \sigma_{r,S}^{k}
\prod_{s \in S} b_{1}^{s}. \label{998}
\end{equation} 
We will refer to $\sigma_{r,S}^{k} (r<0,0 \leq k \leq n-1,S \subseteq
\{1,...,2g\})$ as the Mumford relations.\\
\indent Similarly we know that the restriction of
\[
\Psi^{*}(t) = t^{ng - \delta}c(\pi_{!}({\cal V}^{*} \otimes
\phi^{*}L))(-t^{-1})
\]
to $\HS(\mnd;\Q)$ is a polynomial. As before we may put $\Psi^{*}(t)$ in
the form
\[
\Psi^{*}(t) = \sum_{r=-\infty}^{\g}(\tau_{r}^{0}+\tau_{r}^{1}t+ \cdot
\cdot \cdot + \tau_{r}^{n-1}t^{n-1})(\om(t))^{r}
\]
where $\om(t) = t^{n}+a_{1}t^{n-1}+ \cdot \cdot \cdot +a_{n}$ and
similarly we write
\begin{equation}
\tau_{r}^{k}=\sum_{S \subseteq \{1,...,2g\}} \tau_{r,S}^{k} \prod_{s
\in S} b_{1}^{s}. \label{999}
\end{equation}
We will refer to $\tau_{r,S}^{k} (r<0,0 \leq k \leq n-1,S \subseteq
 \{1,...,2g\})$ as the dual Mumford relations.\\
\indent The motivation for this is that the restrictions of
 $\sigma_{r,S}^{k}$ and $\tau_{r,S}^{k}$ to the strata $\Cmu$ are
 easier to calculate in this form. This is a crucial step in applying
 the following completeness criteria.\\
\indent Given $\mu=(\mu_{1},...,\mu_{n}),\nu=(\nu_{1},...,\nu_{n}) \in
{\cal M}$ then we write $\nu \prec \mu$ if there exists $T$, $1 \leq T
\leq n$, such that
\[
\nu_{i} = \mu_{i} \mbox{ for } T < i \leq n \mbox{ and } \nu_{T} >
\mu_{T}.
\]
We write $\nu \preceq \mu$ if $\nu \prec \mu$ or $\nu = \mu.$ A few
easy calculations verify that $\preceq$ is a total order on ${\cal M}$
with minimal element $\mu_{0}$, the semistable type. For an unstable
type $\mu$ we will write $\mu-1$ for the type previous to $\mu$ with
respect to $\preceq.$
\begin{prop}
\label{KCC}
(Completeness Criteria) Let ${\cal R}$ be a
subset of the kernel of the restriction map
\[
\HG(\C;\Q) \rightarrow \HG(\Css;\Q).
\]
Suppose that for each unstable type $\mu$ there is a subset ${\cal
R}_{\mu}$ of the ideal generated by ${\cal R}$ such that the image of
${\cal R}_{\mu}$ under the restriction map
\[
\HG(\C;\Q) \rightarrow \HG(\C_{\nu};\Q)
\]
is zero when $\nu \prec \mu$ and when $\nu = \mu$ contains the ideal of
$\HG(\C_{\mu};\Q)$ generated by $e_{\mu}$, the equivariant
Euler class of ${\cal N}_{\mu}$, the normal bundle to the stratum $\C_{\mu}$ in
$\C.$ Then ${\cal R}$ generates the kernel of the restriction map
\[
\HG(\C;\Q) \rightarrow \HG(\Css;\Q)
\]
as an ideal of $\HG(\C;\Q).$
\end{prop}
\begin{rem}
The proof of proposition \ref{KCC} below follows similar
lines to the proof of \cite[prop.1]{K2}. However there are some
differences -- the order $\preceq$ does not generally coincide with
$\leq$ -- and further the proof of \cite[p.867]{K2} as given is true
only for the rank two case. For these reasons we include a proof of
proposition \ref{KCC} below although it clearly owes many of its origins to
\cite{K2}.
\end{rem}
{\bf Proof} Let $\mu \in {\cal M}$ and define
\[
V_{\mu} = \bigcup_{\nu \preceq \mu} \C_{\nu}.
\]
We will firstly show that $V_{\mu}$ is an open subset of $\C$
containing $\Cmu$ as a closed submanifold. Note that if $\nu \leq \mu$
then $\nu \preceq \mu$ and thus by property (\ref{11}) if $\nu
\succ \mu$ then $\overline{\C}_{\nu} \subseteq \C - V_{\mu}$. The
stratification is locally finite and hence $V_{\mu}$ is open. Further
note that the closure of $\Cmu$ in $V_{\mu}$ equals
\[
V_{\mu} \cap \bigcup_{\nu \geq \mu} \C_{\nu} = \Cmu
\]
as required.\\
\indent Recall now that the composition of the Thom-Gysin map
\[
H_{\G}^{*-2d_{\mu}}(\C_{\mu};\Q) \rightarrow \HG(V_{\mu};\Q)
\]
with the restriction map
\[
\HG(V_{\mu};\Q) \rightarrow \HG(\C_{\mu};\Q)
\]
is given by multiplication by the Euler class $e_{\mu}$ which is
not a zero-divisor in $\HG(\C_{\mu};\Q)$ \cite[p.569]{AB}.
It follows from the exactness of the Thom-Gysin sequence
\[
\cdot \cdot \cdot \rightarrow H_{\G}^{*-2d_{\mu}}(\C_{\mu};\Q)
\rightarrow \HG(V_{\mu};\Q)
\rightarrow \HG(V_{\mu-1};\Q) \rightarrow \cdot \cdot \cdot
\]
that the direct sum of the restriction maps
\[
\HG(V_{\mu};\Q) \rightarrow \HG(\C_{\mu};\Q) \oplus \HG(V_{\mu-1};\Q)
\]
is injective. Hence inductively the direct sum of restriction maps
\[
\HG(V_{\mu-1};\Q) \rightarrow \bigoplus_{\nu \prec \mu} \HG(\C_{\nu};\Q)
\]
is injective and in particular the image of any element of ${\cal R}_{\mu}$
under the restriction map
\[
\HG(\C;\Q) \rightarrow \HG(V_{\mu-1};\Q)
\]
is zero.\\
\indent For any given $i \geq 0$ there are only finitely many $\nu \in
{\cal M}$ such that $2d_{\nu} \leq i$ and so for each $i \geq 0$ there
exists some $\mu$ such that
\[
H_{\G}^{i}(\C;\Q) = H_{\G}^{i}(V_{\mu};\Q).
\]
Hence it is enough to show that for each $\mu$ the image in
$\HG(V_{\mu};\Q)$ of the ideal generated by ${\cal R}$ contains the
image in $\HG(V_{\mu};\Q)$ of the kernel of the restriction map
\begin{equation}
\HG(\C;\Q) \rightarrow \HG(\Css;\Q). \label{QQQ}
\end{equation}
Note that the above is clearly true for $\mu=\mu_{0}$ as
$V_{\mu_{0}}=\Css.$ We will proceed by induction with respect to $\preceq$.\\
\indent Assume now that $\mu \neq \mu_{0}$ and that $\alpha \in
\HG(\C;\Q)$ lies in the kernel of (\ref{QQQ}).
Suppose that the image of $\alpha$ in $\HG(V_{\mu-1};\Q)$ is in the
image of the ideal generated by ${\cal R}.$ We may, without any loss of
generality, assume that the image of $\alpha$ in $\HG(V_{\mu-1};\Q)$ is
zero. Thus by the exactness of the Thom-Gysin sequence there exists an
element $\beta \in H_{\G}^{*-2d_{\mu}}(\C_{\mu};\Q)$ which is
mapped to the image of $\alpha$ in $\HG(V_{\mu};\Q)$ by the Thom-Gysin
map
\[
H_{\G}^{*-2d_{\mu}}(\C_{\mu};\Q) \rightarrow \HG(V_{\mu};\Q).
\]
Hence the image of $\alpha$ under the restriction map
\[
\HG(\C;\Q) \rightarrow \HG(\C_{\mu};\Q)
\]
is $\beta e_{\mu}$, and by hypothesis there is an element $\gamma$ of ${\cal
R}_{\mu}$ which maps under the
restriction map
\[
\HG(V_{\mu};\Q) \rightarrow \HG(\C_{\mu};\Q)
\]
to $\beta e_{\mu}.$ Now the images of $\gamma$ and $\alpha$ in
$\HG(V_{\mu-1};\Q)$ are both zero and we also know the direct sum of the
restriction maps
\[
\HG(V_{\mu};\Q) \rightarrow \HG(\C_{\mu};\Q) \oplus \HG(V_{\mu-1};\Q)
\]
to be injective. Thus the images of $\gamma$ and $\alpha$ in
$\HG(V_{\mu};\Q)$ are the same, completing the proof. $\indent \Box$
\begin{rem}
Kirwan's completeness criteria follow from the
above criteria since for each $\mu$
\[
V_{\mu-1} \subseteq \C - \bigcup_{\nu \geq \mu} \Cnu.
\]
So if the restriction of a relation to $\HG(\Cnu;\Q)$ vanishes for
every $\nu \not \geq \mu$ then certainly the same relation restricts
to zero in $\HG(\Cnu;\Q)$ for any $\nu \prec \mu.$
\end{rem}
\begin{rem}
Kirwan's proof of Mumford's
conjecture \cite[$\S$ 2]{K2} amounts to showing that for each unstable
type $\mu=(d_{1},d_{2})$ the set
\[
{\cal R}_{\mu} = \bigcup \{
\sigma_{d_{2}-2g+1,S}^{0},\sigma_{d_{2}-2g+1,S}^{1} \},
\]
where the union is taken over all subsets $S \subseteq \{1,...,2g\}$,
satisfies the above criteria. In the rank two case the criteria of
proposition \ref{KCC} are in fact equivalent to Kirwan's completeness
criteria since $\preceq$ and $\leq$ coincide.
\end{rem}
\section{Chern Class Computations.}
\indent We first describe the restriction maps $\HG(\C;\Q) \rightarrow
\HG(\Cmu;\Q)$ and our preferred generators for
$\HG(\Cmu;\Q)$. Let $\mu=(d_{1}/n_{1},...,d_{P}/n_{P})$. Let
$\C(n_{p},d_{p})^{ss}$ denote the space of all semistable holomorphic
structures on a fixed Hermitian vector bundle of rank $n_{p}$ and
degree $d_{p}$ and let $\G(n_{p},d_{p})$ be the gauge group of that
bundle. Atiyah and
Bott \cite[prop. 7.12]{AB} show that the map
\[
\prod_{p=1}^{P} \C(n_{p},d_{p})^{ss} \rightarrow \Cmu,
\]
which sends a sequence of semistable bundles $(F_{1},...,F_{P})$ to
the direct sum $F_{1} \oplus \cdot \cdot \cdot \oplus F_{P}$, induces an
isomorphism
\[
\HG(\Cmu;{\bf Q}) \cong \bigotimes_{1 \leq p \leq P} H^{*}_{{\cal
G}(n_{p},d_{p})}({\cal C}(n_{p},d_{p})^{ss};{\bf Q}).
\]
Thus we can find generators
\begin{equation}
\bigcup_{p=1}^{P} \left( \{a_{r}^{p} | 1 \leq r \leq n_{p} \}
\cup \{ b_{r}^{p,s} | 1 \leq r \leq n_{p},1 \leq s \leq 2g \} \cup \{
f_{r}^{p} | 2 \leq r \leq n_{p} \} \right) \label{A}
\end{equation}
corresponding to the generators of $\HG(\Css;\Q)$ described earlier in
(\ref{0}). As before we also define
\[
\xi_{i,j}^{p,q} = \sum_{s=1}^{g} b_{i}^{p,s} b_{j}^{q, s+g}.
\]
To explicitly describe the restriction map note that
$c_{r}({\cal V})$ restricts to $c_{r}(\bigoplus_{p=1}^{P} {\cal
V}_{p})$ where ${\cal V}_{p}$ is the universal bundle on $\C(n_{p},d_{p})$. The
restrictions of the generators of $\HG(\C;\Q)$ can be written in terms of
the generators of $\HG(\Cmu;\Q)$ by taking the appropriate
coefficients in the K\"{u}nneth decomposition.\\
\indent One problem that we will be faced with in due course is how to
calculate the coefficients of $\prod_{s \in S} b_{1}^{s}$ once we have
restricted to a stratum. Suppose first that the stratum concerned is of type
$\mu=(d_{1},...,d_{n}) \in \Delta$ and take
$\zeta \in \HG(\C;\Q).$ We can express $\zeta$ in terms of the generators
\[
\{a_{r} | 1 \leq r \leq n \} \cup \{ b_{r}^{s} | 1 \leq r \leq n,1
\leq s \leq 2g \} \cup \{ f_{r} | 2 \leq r \leq n \}
\]
but equally we could write $\zeta$ in terms of
\[
\{a_{r} | 1 \leq r \leq n \} \cup \{ n b_{r}^{s}
-(n-r+1)a_{r-1}b_{1}^{s} | 2 \leq r \leq n,1
\leq s \leq 2g \}
\]
\begin{equation}
\cup \{ n^{2} f_{r} - n(n-r+1)(\xi_{r-1,1} + \xi_{1,r-1}) + (n-r+1)(n-r+2) a_{r-2} \xi_{1,1} | 2 \leq r
\leq n \} \label{1000}
\end{equation}
and $\{b_{1}^{s} | 1 \leq s \leq 2g\}.$ We shall take the coefficients
of $\prod_{s \in S}b_{1}^{s}$ when $\zeta$ is expressed in this latter
form. The reason for this is that
the restrictions of the elements (\ref{1000}) in $\HG(\Cmu;\Q)$ can then
be written in terms of
\begin{equation}
\{a_{1}^{r} | 1 \leq r \leq n \} \cup \{b_{1}^{p,s}-b_{1}^{n,s} | 1
\leq p \leq n-1, 1 \leq s \leq 2g \} \label{1001}
\end{equation}
(see remark \ref{tedious}.) We can uniquely write the restriction of $\zeta$ in terms of the
elements (\ref{1001}) and the restrictions of $b_{1}^{s},(1 \leq s
\leq 2g)$. Hence we may calculate the restrictions of the coefficients of
$\prod_{s \in S} b_{1}^{s}$ in $\zeta$ by taking the coefficients of
\[
\prod_{s \in S} (b_{1}^{1,s} + \cdots + b_{1}^{n,s})
\]
in the restriction of $\zeta$.\\
\indent We deal with a general type stratum in a similar way. Let
$\mu = (d_{1}/n_{1},...,d_{P}/n_{P})$. We define formal symbols
$a^{p,k},b^{p,k,s}$ and $d^{p,k}$ such that the $r$th Chern class
$c_{r}({\cal V}_{p})$ is given by the $r$th elementary symmetric
polynomial in
\begin{equation}
a^{p,k} + \sum_{s=1}^{2g} b^{p,k,s} \otimes \alpha_{s} + d^{p,k}
\otimes \omega \indent (1 \leq k \leq n_{p}) \label{QQ}
\end{equation}
when $1 \leq r \leq n_{p}$ and $1 \leq p \leq P$. In terms of
$a^{p,k},b^{p,k,s}$ and $d^{p,k}$ the restriction map to
$\HG(\Cmu;\Q)$ is formally the same as the restriction map when $\mu \in
\Delta$. Again we may uniquely write the restriction of
$\zeta$ in terms of
\begin{equation}
\bigcup_{p=1}^{P} \bigcup_{k=1}^{n_{p}} \{a^{p,k}, d^{p,k} \}
\cup \bigcup_{p=1}^{P-1} \bigcup_{k= 1}^{n_{p}} \bigcup_{s=1}^{2g} \{
b^{p,k,s} - b^{P,n_{P},s} \} \cup \bigcup_{k=
1}^{n_{P}-1} \bigcup_{s=1}^{2g} \{ b^{P,k,s} - b^{P,n_{P},s} \} \label{NEW1}
\end{equation}
and the restrictions of $b_{1}^{s}, (1 \leq s \leq 2g)$, and we take
the coefficients of
\[
\prod_{s \in S}(b_{1}^{1,s} + \cdots + b_{1}^{P,s})
\]
as before.\\
\indent So in our definitions of the Mumford and dual Mumford
relations, (\ref{998}) and (\ref{999}), we assume first that
$\sigma_{r}^{k}$ and $\tau_{r}^{k}$ have first been written in terms
of the elements (\ref{1000}) before taking the appropriate
coefficient.
\begin{rem}
\label{tedious}
It is a trivial but tedious calculation to show that the restrictions
of the elements (\ref{1000}) in $\HG(\Cmu;\Q)$ for $\mu \in \Delta$ can
indeed be written in terms of the elements (\ref{1001}). Let
$a_{r}^{\mu}$ denote the restriction of $a_{r}$ to $\HG(\Cmu;\Q)$;
this equals the $r$th elementary symmetric product in
$a_{1}^{1},...,a_{1}^{n}$. The restrictions of $b_{r}^{s}$ and $f_{r}$
in $\HG(\Cmu;\Q)$ equal
\[
\sum_{i=1}^{n} b_{1}^{i,s} \frac{\partial a_{r}^{\mu}}{\partial
a_{1}^{i}}, \indent \sum_{i=1}^{n} d_{i} \frac{\partial
a_{r}^{\mu}}{\partial a_{1}^{i}} + \sum_{i=1}^{n} \sum_{j=1}^{n}
\xi_{1,1}^{i,j} \frac{\partial^{2} a_{r}^{\mu}}{\partial a_{1}^{i}
\partial a_{1}^{j}}.
\]
The restrictions of the elements (\ref{1000}) can then be seen to
equal
\[
a_{r}^{\mu}, \indent \sum_{i=1}^{n-1} \left( n \frac{\partial
a_{r}^{\mu}}{\partial a_{1}^{i}} - (n-r+1) a_{r-1}^{\mu} \right) (
b_{1}^{i,s} - b_{1}^{n,s}),
\]
and 
\[
n^{2} \sum_{i=1}^{n} d_{i} \frac{\partial a_{r}^{\mu}}{\partial
a_{1}^{i}} + \sum_{i=1}^{n-1} \sum_{j=1}^{n-1} \sum_{s=1}^{g}
(b_{1}^{i,s} - b_{1}^{n,s})(b_{1}^{j,s+g} - b_{1}^{n,s+g}) \left(
n^{2} \frac{\partial^{2}a_{r}^{\mu}}{\partial a_{1}^{i} \partial
a_{1}^{j}} \right.
\]
\[
\left. - n(n-r+1) \left(\frac{\partial a_{r-1}^{\mu}}{\partial
a_{1}^{i}} + \frac{\partial a_{r-1}^{\mu}}{\partial a_{1}^{j}} \right) +
(n-r+1)(n-r+2)a_{r-2}^{\mu} \right).
\]  
\end{rem}
\indent The remains of this section are given over to calculating the
Mumford and dual Mumford relations. Our first problem is to obtain
their generating functions from their respective Chern characters which we can
evaluate using the Grothendieck-Riemann-Roch theorem (GRR).
\begin{lem}
\label{Chernlemma}
Suppose that
\begin{equation}
{\rm ch}(E) = \sum_{i=1}^{M} \alpha_{i} e^{\delta_{i}} + \sum_{i=1}^{N}
\beta_{i} e^{\epsilon_{i}} \label{230}
\end{equation}
where the $\beta_{i},\delta_{i}$ and the $\epsilon_{i}$ are formal
degree two classes and the $\alpha_{i}$ are formal degree zero
classes. Then as a formal power series
\begin{equation}
c(E)(t) = \sum_{r=0}^{\infty} c_{r}(E) \cdot t^{r} = \prod_{i=1}^{M}
(1+\delta_{i}t)^{\alpha_{i}} \prod_{i=1}^{N} \exp \left\{
\frac{\beta_{i}t}{1+\epsilon_{i}t} \right\}. \label{231}
\end{equation}
\end{lem}
{\bf Proof}
The relationship between the Chern character and Chern polynomial
is as follows. If $\ch(E)= \sum_{i=1}^{K} e^{\gamma_{i}}$ where $\gamma_{i}$
are formal degree two classes then
\[
c(E)(t) = \prod_{i=1}^{K} (1+\gamma_{i}t).
\]
If $\ch(E)$ is in the form of (\ref{230}) then by comparing
degrees we find that
\[
\sum_{i=1}^{M} \alpha_{i}(\delta_{i})^{n} + \sum_{i=1}^{N} n \beta_{i}
(\epsilon_{i})^{n-1} = \sum_{i=1}^{K} (\gamma_{i})^{n}
\]
for each $n \geq 0$. Thus on
the level of formal power series $\log c(E)(t)$ equals
\[
\sum_{i=1}^{K}  \sum_{r=1}^{\infty} (-1)^{r+1}
\frac{(\gamma_{i}t)^{r}}{r} = \sum_{i=1}^{M} \alpha_{i} \log(1+\delta_{i}t) +
\sum_{i=1}^{N} \frac{\beta_{i}t}{1+\epsilon_{i}t}
\]
and hence the result (\ref{231}). $\indent \Box$.\\[\baselineskip]
\indent Armed with the above lemma we are now in a position to
determine the Chern polynomials $c(\pi_{!}{\cal V})(t)$ and
$c(\pi_{!}({\cal V}^{*} \otimes \phi^{*}L))(-t)$. We can, and
will, calculate these Chern polynomials in terms of the generators
$a_{r},b_{r}^{s}$ and $f_{r}$ of $\HG({\cal C};\Q)$ (see (\ref{26})
and (\ref{210})). However the expressions obtained are
somewhat cumbersome and for ease of calculation we
will find the formal expressions, (\ref{25}) and (\ref{29}),
calculated directly from the above lemma of more use.
\begin{prop}
\label{Chernprop}
The Chern polynomial $c(\pi_{!}{\cal V})(t)$ equals
\begin{equation}
\Omega(t)^{-\g} \prod_{k=1}^{n} (1+\delta_{k}t)^{W_{k}} \exp \left\{
\frac{X_{k}t}{1+\delta_{k}t} \right\} \label{25}
\end{equation}
and $c(\pi_{!}({\cal
V}^{*} \otimes \phi^{*}L))(-t)$ equals
\begin{equation}
\Omega(t)^{3\g+1} \prod_{k=1}^{n} (1+\delta_{k}t)^{-W_{k}} \exp
\left\{ \frac{-X_{k}t}{1+\delta_{k}t} \right\}, \label{29}
\end{equation}
where $\delta_{1},...,\delta_{n}$ are formal degree two classes such
that their $r$th elementary symmetric polynomial equals $a_{r},$ and
\[
\Omega(t)= \prod_{k=1}^{n}(1+\delta_{k}t) = 1 + a_{1}t + \cdot \cdot
\cdot + a_{n}t^{n}, \indent \xi_{i,j}= \sum_{s=1}^{g} b_{i}^{s}
b_{j}^{s+g},
\]
\[
W_{k}= \sum_{i=1}^{n} f_{i} \frac{\partial \delta_{k}}{\partial a_{i}} +
\sum_{i=1}^{n} \sum_{j=1}^{n} \xi_{i,j} \frac{\partial^{2} \delta_{k}}{\partial
a_{i} \partial a_{j}}, \indent X_{k} = \sum_{i=1}^{n} \sum_{j=1}^{n} \xi_{i,j}
\frac{\partial \delta_{
k}}{\partial a_{i}} \frac{\partial \delta_{k}}{\partial a_{j}}.
\]
In terms of the generators $a_{r},b_{r}^{s}$ and $f_{r}$ for
$\HG({\cal C};{\bf Q})$ then $c(\pi_{!}{\cal V})(t)$ equals
\begin{equation}
\Omega(t)^{-\g} \exp \left\{ \int_{0}^{t} \left( \frac{d}{u} - \sum_{i=1}^{n}
\frac{f_{i} u^{i-2}}{\Omega(u)} + \sum_{i=1}^{n} \sum_{j=1}^{n}
\frac{\xi_{i,j}u^{i+j-2}}{\Omega(u)^{2}} \right) {\rm d}u \right\} \label{26}
\end{equation}
and $c(\pi_{!}({\cal V}^{*} \otimes \phi^{*}L))(-t)$ equals
\begin{equation}
\Omega(t)^{3\g+1} \exp  \left\{ \int_{0}^{t} \left( -\frac{d}{u} + \sum_{i=1}^{n}
\frac{f_{i} u^{i-2}}{\Omega(u)} - \sum_{i=1}^{n} \sum_{j=1}^{n}
\frac{\xi_{i,j}u^{i+j-2}}{\Omega(u)^{2}} \right) {\rm d}u \right\}. \label{210}
\end{equation}
\end{prop}
{\bf Proof}
Now $\ch({\cal V}) = e^{\gamma_{1}} + \cdot \cdot \cdot +
e^{\gamma_{n}}$ where $\gamma_{1},...,\gamma_{n}$ are formal degree two classes
such that their $r$th elementary symmetric polynomial equals
\[
c_{r}({\cal V}) = a_{r} \otimes 1 + \sum_{s=1}^{2g} b_{r}^{s} \otimes
\alpha_{s} + f_{r} \otimes \omega \indent (1 \leq r \leq n) .
\]
For each $k \geq 0$ there exist coefficients $\rho_{r_{1},...,r_{n}}^{(k)}$
such that
\[
(\gamma_{1})^{k} + \cdot \cdot \cdot + (\gamma_{n})^{k} = \sum
\rho_{r_{1},...,r_{n}}^{(k)} (c_{1}({\cal V}))^{r_{1}} \cdot \cdot \cdot
(c_{n}({\cal V}))^{r_{n}}
\]
where the sum is taken over all non-negative $r_{1},...,r_{n}$ such that
$r_{1}+2r_{2}+ \cdots +nr_{n} = k$. Now
\[
(a_{1} \otimes 1 + \sum_{s=1}^{2g} b_{1}^{s} \otimes \alpha_{s} + f_{1} \otimes
\omega)^{r_{1}} \cdot \cdot \cdot (a_{n} \otimes 1 + \sum_{s=1}^{2g} b_{n}^{s}
\otimes \alpha_{s} + f_{n} \otimes \omega)^{r_{n}}
\]
equals
\[
(a_{1})^{r_{1}} \cdot \cdot \cdot (a_{n})^{r_{n}} \otimes 1 + \sum_{i=1}^{n}
\sum_{s=1}^{2g} b_{i}^{s} \dai (a_{1})^{r_{1}} \cdot \cdot \cdot
(a_{n})^{r_{n}} \otimes \alpha_{s}
\]
\[
+\sum_{i=1}^{n} f_{i} \dai (a_{1})^{r_{1}} \cdot \cdot \cdot (a_{n})^{r_{n}}
\otimes \omega + \sum_{i=1}^{n} \sum_{j=1}^{n} \xi_{i,j} \daij (a_{1})^{r_{1}}
\cdot \cdot \cdot (a_{n})^{r_{n}} \otimes \omega.
\]
Since
\[
\sum \rho_{r_{1},...,r_{n}}^{(k)} (a_{1})^{r_{1}} \cdot \cdot \cdot
(a_{n})^{r_{n}} = (\delta_{1})^{k} + \cdot \cdot \cdot + (\delta_{n})^{k}
\]
we find that $\ch({\cal V})$ equals
\[
\sum_{k=1}^{n} e^{\delta_{k}} \otimes 1 + \sum_{i=1}^{n} \sum_{s=1}^{2g}
\sum_{k=1}^{n} b_{i}^{s} \dai e^{\delta_{k}} \otimes \alpha_{s}
\]
\begin{equation}
+\sum_{i=1}^{n} \sum_{k=1}^{n} f_{i} \dai e^{\delta_{k}} \otimes
\omega + \sum_{i=1}^{n} \sum_{j=1}^{n} \sum_{k=1}^{n} \xi_{i,j} \daij
e^{\delta_{k}} \otimes \omega. \label{212}
\end{equation}
{}From GRR we have $\ch(\pi_{!}{\cal V}) = \pi_{*} ( \ch({\cal V}) \cdot
1 \otimes (1-\g \omega))$ and hence $\ch(\pi_{!}{\cal V})$ equals
\[
\sum_{i=1}^{n} \sum_{k=1}^{n} f_{i} \dai e^{\delta_{k}} + \sum_{i=1}^{n}
\sum_{j=1}^{n} \sum_{k=1}^{n} \xi_{i,j} \daij e^{\delta_{k}} - \g
\sum_{k=1}^{n} e^{\delta_{k}}
= \sum_{k=1}^{n} (-\g + W_{k} + X_{k}) e^{\delta_{k}}.
\]
Note that $W_{k}$ has degree zero and $X_{k}$ has degree two. Hence by
lemma \ref{Chernlemma} we see that $c(\pi_{!}{\cal V})(t)$ equals
\[
(\Omega(t))^{-\g} \prod_{k=1}^{n} (1+\delta_{k}t)^{W_{k}} \exp \left\{
\frac{X_{k}t}{1+\delta_{k}t} \right\}
\]
to give equation (\ref{25}).\\
\indent Now $\frac{{\rm d}}{{\rm d}t} \log ( \Omega(t)^{\g}
c(\pi_{!}{\cal V})(t))$ equals
\[
\sum_{i=1}^{n} \sum_{k=1}^{n} f_{i} \frac{\partial \delta_{k}}{\partial a_{i}}
\frac{\delta_{k}}{1+\delta_{k}t} + \sum_{i=1}^{n} \sum_{j=1}^{n} \sum_{k=1}^{n}
\xi_{i,j} \left( \frac{\partial^{2} \delta_{k}}{\partial a_{i} \partial a_{j}}
\frac{\delta_{k}}{1+\delta_{k} t} + \frac{\partial \delta_{k}}{\partial a_{i}}
\frac{\partial \delta_{k}}{\partial a_{j}} \frac{1}{(1+ \delta_{k} t)^{2}}
\right)
\]
\[
= \sum_{i=1}^{n} \sum_{j=1}^{n} \frac{\xi_{i,j}}{t^{2}} \left( \sum_{k=1}^{n} t
\frac{\partial^{2} \delta_{k}}{\partial a_{i} \partial a_{j}} - \sum_{k=1}^{n}
\left( \frac{t}{1+ \delta_{k}t} \frac{\partial^{2} \delta_{k}}{\partial a_{i}
\partial a_{j}} -
\frac{t^{2}}{(1+\delta_{k} t)^{2}} \frac{\partial \delta_{k}}{\partial a_{i}}
\frac{\partial \delta_{k}}{\partial a_{j}} \right) \right)
\]
\begin{equation}
+ \sum_{i=1}^{n} \frac{f_{i}}{t} \left( \sum_{k=1}^{n} \frac{\partial
\delta_{k}}{\partial a_{i}} - \sum_{k=1}^{n} \frac{\partial
\delta_{k}}{\partial a_{i}} \frac{1}{1+ \delta_{k} t} \right) \label{219}
\end{equation}
Since $\sum_{k=1}^{n} \frac{\partial \delta_{k}}{\partial a_{i}} =
\frac{ \partial a_{1}}{\partial a_{i}}$, $f_{1}=d$, and
$\sum_{k=1}^{n} \frac{\partial^{2} \delta_{k}}{\partial a_{i} \partial
a_{j}} = \frac{\partial^{2} a_{1}}{\partial a_{i} \partial a_{j}} = 0$
then (\ref{219}) reduces to
\[
\frac{d}{t} - \sum_{i=1}^{n} \sum_{k=1}^{n} f_{i} \dai \frac{\log (1+
\delta_{k} t)}{t^{2}} - \sum_{i=1}^{n} \sum_{j=1}^{n} \sum_{k=1}^{n} \xi_{i,j}
\daij \frac{\log (1 + \delta_{k}t)}{t^{2}}
\]
\[
= \frac{d}{t} - \left(
\sum_{i=1}^{n} f_{i} \dai + \sum_{i=1}^{n} \sum_{j=1}^{n} \xi_{i,j}
\daij \right) \frac{\log \Omega(t)}{t^{2}}
\]
to give equality (\ref{26}).\\
\indent The calculations for the dual case follow in a similar
fashion. We have that  $\ch({\cal V}^{*}) = e^{-\gamma_{1}} + \cdot
\cdot \cdot + e^{-\gamma_{n}}$ with $\gamma_{1},...,\gamma_{n}$ as
before and arguing as in the calculation of (\ref{212}) we determine that
$\ch({\cal V}^{*})$ equals
\[
\sum_{k=1}^{n} e^{-\delta_{k}} \otimes 1 + \sum_{i=1}^{n} \sum_{s=1}^{2g}
\sum_{k=1}^{n} b_{i}^{s} \dai e^{-\delta_{k}} \otimes \alpha_{s}
\]
\begin{equation}
+\sum_{i=1}^{n} \sum_{k=1}^{n} f_{i} \dai e^{-\delta_{k}} \otimes
\omega + \sum_{i=1}^{n} \sum_{j=1}^{n} \sum_{k=1}^{n} \xi_{i,j} \daij
e^{-\delta_{k}} \otimes \omega. \label{213}
\end{equation}
We know that $\ch(\phi^{*}L) = \phi^{*}(e^{(4\g+1)\omega}) = 1
\otimes (1+ (4\g+1)\omega)$ and GRR shows that $\ch(\pi_{!}({\cal
V}^{*} \otimes \phi^{*}L))$ equals
\[ \pi_{*}(\ch({\cal
V}^{*}) \cdot \ch(\phi^{*}L) \cdot 1 \otimes (1 - \g \omega)) =
\pi_{*}(\ch({\cal
V}^{*}) \cdot 1 \otimes (1 + (3\g+1) \omega))
\]
which gives
\begin{equation}
\ch(\pi_{!}({\cal V}^{*} \otimes \phi^{*}L)) = \sum_{k=1}^{n} (
(3\g+1) -W_{k} +X_{k} ) e^{-\delta_{k}}. \label{215}
\end{equation}
Applying lemma \ref{Chernlemma} to expression (\ref{215}) gives equation
(\ref{29}). Expression (\ref{210}) is arrived at by calculating $\frac{{\rm
d}}{{\rm d}t} \log ((\Omega(t))^{-3\g-1}c(\pi_{!}({\cal V}^{*} \otimes
\phi^{*}L))(t))$ and grouping the terms in a similar manner to
expression  (\ref{219}). $\indent \Box$
\begin{rem}
Note that $\delta_{k}, W_{k}$ and $X_{k}$ are not
elements of $\HG(\C;\Q)$. However the direct sum of the restriction
maps
\[
\HG(\C;\Q) \rightarrow \bigoplus_{\mu \in \Delta} \HG(\Cmu;\Q)
\]
is injective and so we may consider $\delta_{k},W_{k}$ and $X_{k}$ as
elements of $\bigoplus_{\mu \in \Delta} \HG(\Cmu;\Q)$ corresponding
respectively to $a_{1}^{k},d_{k}$ and $\xi_{1,1}^{k,k}$ in each
summand $\HG(\Cmu;\Q).$
\end{rem}
\begin{rem}
From (\ref{26}) we can find an expression for
\[
\frac{\Psi'(t)}{\Psi(t)} = \frac{d-n\g}{t} - \frac{c(\pi_{!}{\cal
V})'(t^{-1})}{t^{2} c(\pi_{!}{\cal V})(t^{-1})}.
\]
In fact we may write $\Psi'(t)/\Psi(t)$ as a rational function with
denominator $(\om(t))^{2}$ and a numerator of degree at most $2n-1$.
By multiplying by $\Psi(t)$ and comparing coefficients of
$t^{k} (\om(t))^{r}, (r \leq \g, 0 \leq k <n)$ we may derive
recurrence relations amongst the Mumford relations which determine
$\{\sigma_{r}^{k} : 0 \leq k < n\}$ in terms of
$\{\sigma_{r+1}^{k},\sigma_{r+2}^{k} : 0 \leq k < n \}$. Similar
recurrence relations exist among the dual Mumford relations which
determine $\{\tau_{r}^{k} : 0 \leq k < n\}$ in terms of
$\{\tau_{r+1}^{k},\tau_{r+2}^{k} : 0 \leq k < n \}$.
\end{rem}
\indent The calculation of the restriction of $c(\pi_{!}{\cal V})(t)$ to
$\HG(\Cmu;\Q)[[t]]$ follows easily from the previous proposition. As in
\cite[prop. 2]{K2} this restriction can be expressed in terms of
elementary functions of the generators of $\HG(\Cmu;\Q)$ when $\mu \in \Delta$.
However for a general type $\mu$
this restriction cannot be expressed so easily and we will find
formal expressions similar to (\ref{25}) of more use.
\begin{cor}
\label{Chernresn}
Let $\mu= (d_{1}/n_{1},...,d_{P}/n_{P})$. The restriction to $H_{\cal G}^{*}({\cal C}_{\mu};{\bf Q})[[t]]$ of
$c(\pi_{!}{\cal V})(t)$ equals the formal power series
\begin{equation}
\Omega_{\mu}(t)^{-\g} \prod_{p=1}^{P} \prod_{k=1}^{n_{p}}
(1+\delta_{k}^{p}t)^{W_{k}^{p}} \exp \left\{
\frac{X_{k}^{p}t}{1+\delta_{k}^{p}t} \right\} \label{27}
\end{equation}
and similarly the restriction of $c(\pi_{!}({\cal V}^{*} \otimes
\phi^{*}L))(-t)$ to $H_{\cal G}^{*}({\cal C}_{\mu};{\bf Q})[[t]]$ equals
\begin{equation}
\Omega_{\mu}(t)^{3\g+1} \prod_{p=1}^{P} \prod_{k=1}^{n_{p}}
(1+\delta_{k}^{p}t)^{-W_{k}^{p}} \exp \left\{
\frac{-X_{k}^{p}t}{1+\delta_{k}^{p}t} \right\}, \label{211}
\end{equation}
where $\delta_{1}^{p},...,\delta_{n_{p}}^{p}$ are formal degree two
classes such that their $r$th elementary symmetric polynomial equals
$a_{r}^{p}$, where $\Omega_{\mu}(t)= \prod_{p=1}^{P} \prod_{k=1}^{n_{p}} (1+
\delta_{k}^{p}t)$ is the restriction of $\Omega(t)$ to $\HG(\Cmu;{\bf
Q})[t],$ and where $\xi_{i,j}^{p,p},W_{k}^{p}$ and $X_{k}^{p}$
correspond to the expressions defined in the statement of proposition
\ref{Chernprop}.
\end{cor}
{\bf Proof}
Expression (\ref{27}) is immediate from the previous proposition
once we note that the restriction of $\ch(\pi_{!}{\cal V})$ to $H^{*}_{\cal G}
( {\cal
C}_{\mu};{\bf Q})$ equals
\[
\sum_{p=1}^{P} \pi_{*}(\ch({\cal V}_{p}) \cdot 1 \otimes (1-\g \omega))
\]
and recall that the Chern polynomial is
multiplicative. The dual expression (\ref{211}) follows in a similar
fashion. $\Box$
\begin{cor}
\label{Deltaresn}
Let $\mu=(d_{1},...,d_{n}) \in \Delta$. Then the restriction
of $c(\pi_{!}{\cal V})(t)$ to
$H^{*}_{\cal G} ({\cal C}_{\mu};{\bf Q})[[t]]$ equals
\[
\prod_{p=1}^{n} (1+ a_{1}^{p} t)^{d_{p} -\g} \exp \left\{ \frac{\xi_{1,1}^{p,p}
t}{1 + a_{1}^{p} t} \right\}.
\]
Also the restriction of $c(\pi_{!}({\cal V}^{*}
\otimes \phi^{*}L))(-t)$ to
$H^{*}_{\cal G} ({\cal C}_{\mu};{\bf Q})[[t]]$ equals
\[
\prod_{p=1}^{n} (1+ a_{1}^{p} t)^{3\g+1-d_{p}} \exp \left\{
\frac{-\xi_{1,1}^{p,p} t}{1 + a_{1}^{p} t} \right\}.
\]
\end{cor}
{\bf Proof}
Simply note that in this case $\delta_{1}^{p} =
a_{1}^{p},W_{1}^{p} = d_{p}$ and $X_{1}^{p}=\xi_{1,1}^{p,p}. \indent
\Box$
\begin{rem}
Let $\mu=(d_{1}/n_{1},...,d_{P}/n_{P})$. From the calculation
(\ref{212}) and since the Chern
character is additive we know that the restriction of $\ch({\cal V})$
to $\HG(\Cmu;\Q)$ equals
\[
\sum_{p=1}^{P} \sum_{k=1}^{n_{p}} \exp \left\{ \delta_{k}^{p} +
\sum_{s=1}^{2g} \left( \sum_{i=1}^{n_{p}} b_{i}^{p,s}\frac{\partial
\delta_{k}^{p}}{\partial a_{i}^{p}} \right) \otimes \alpha_{s} + W_{k}^{p}
\otimes \omega \right\}.
\]
Thus in terms of our earlier notation (\ref{QQ}) we have
\[
a^{p,k} = \delta_{k}^{p}, \indent b^{p,k,s} = \sum_{i=1}^{n_{p}}
b_{i}^{p,s} \frac{\partial
\delta_{k}^{p}}{\partial a_{i}^{p}}, \indent d^{p,k}=W_{k}^{p}.
\]
\end{rem}
\indent We end this section with two further calculations, namely the Chern
polynomials of the normal bundle ${\cal N}_{\mu}$ to the stratum $\Cmu$ in $\C$
(necessary to the completeness criteria) and of the tangent bundle $T$
to the moduli space $\mnd$ (needed for generalising the proof of the Newstead-Ramanan conjecture).
\begin{lem}
\label{normal}
Let $\mu=(d_{1}/n_{1},...,d_{P}/n_{P})$. Then the
Chern polynomial $c({\cal N}_{\mu})(t)$ of the normal bundle in ${\cal
C}$ to the stratum $\Cmu$ equals
\begin{equation}
{\cal P}_{\mu}(t)^{\g} \prod_{I<J} \prod_{k=1}^{n_{I}}
\prod_{l=1}^{n_{J}}
(1+(\delta_{l}^{J}-\delta_{k}^{I})t)^{W_{k}^{I}-W_{l}^{J}} \exp
\left\{ \frac{-\Xi_{k,l}^{I,J}t}{1+(\delta_{l}^{J}-\delta_{k}^{I})t}
\right\} \label{216}
\end{equation}
where
\[
\Xi_{k,l}^{I,J}= \sum_{s=1}^{g} \left( \sum_{i=1}^{n_{I}} b_{i}^{I,s}
\frac{\partial \delta_{k}^{I}}{\partial a_{i}^{I}} -\sum_{j=1}^{n_{J}}
b_{j}^{J,s} \frac{\partial \delta_{l}^{J}}{\partial a_{j}^{J}} \right) \left(
\sum_{i=1}^{n_{I}} b_{i}^{I,s+g} \frac{\partial \delta_{k}^{I}}{\partial a_{i}^{I}} -\sum_{j=1}^{n_{J}}
b_{j}^{J,s+g} \frac{\partial \delta_{l}^{J}}{\partial a_{j}^{J}} \right)
\]
and
\[
{\cal P}_{\mu}(t)= \prod_{I<J} \prod_{k=1}^{n_{I}} \prod_{l=1}^{n_{J}}
(1+(\delta_{l}^{J} - \delta_{k}^{I})t).
\]
\end{lem}
{\bf Proof}
Kirwan \cite[lemma 2]{K2} showed that the normal bundle
${\cal N}_{\mu}$ to $\Cmu$ in $\C$, equals
\[
-\pi_{!} \left( \bigoplus_{I < J} {\cal V}^{*}_{I} \otimes {\cal
V}_{J} \right).
\]
{}From the proof of the proposition \ref{Chernprop} we can find expressions for
$\ch({\cal V}_{J})$ and $\ch({\cal V}_{I}^{*})$ corresponding to (\ref{212})
and
(\ref{213}). The GRR implies that
\[
\ch({\cal N}_{\mu}) = \sum_{I<J} \pi_{*} ( \ch({\cal V}_{I}^{*}) \cdot
\ch({\cal V}_{J}) \cdot 1 \otimes (\g \omega -1)).
\]
Substituting in these expressions for $\ch({\cal V}_{J})$ and $\ch({\cal
V}_{I}^{*})$ we find that $\ch({\cal N}_{\mu})$ equals
\[
\sum_{I<J} \left\{ \sum_{k=1}^{n_{I}} \sum_{l=1}^{n_{J}}
(\g+W_{k}^{I}-W_{l}^{J}-\Xi_{k,l}^{I,J})e^{\delta_{l}^{J}-\delta_{k}^{I}}
\right\}.
\]
Applying lemma \ref{Chernlemma} produces the required result (\ref{216}). $\indent
\Box$
\begin{lem}
\label{Pont}
The total Pontryagin class of $\mnd$ equals
\[
\prod_{1 \leq k < l \leq n} (1 + (\delta_{k} - \delta_{l})^{2})^{2\g}.
\]
In particular the Pontryagin ring of $\mnd$ is generated by the elementary
symmetric polynomials in
\[
\{ ( \delta_{k} -\delta_{l}) ^{2} : 1 \leq k < l \leq n \}.
\]
\end{lem}
{\bf Proof}
Let $T$ denote the tangent bundle of $\mnd$. From \cite[p.582]{AB} we
know that
\[
T + T^{*} -2 = \pi_{!}({\rm End} V \otimes (\Omega_{M}^{1}-1)).
\]
Applying GRR we find
\[
\ch T + \ch T^{*} - 2 = 2 \g \ch ({\rm End V}|\mnd)
\]
which we know to equal
\[
2 \g \left(\sum_{k=1}^{n} e^{\delta_{k}} \right) \left( \sum_{l=1}^{n}
e^{-\delta_{l}} \right)
\]
from expressions (\ref{212}) and (\ref{213}).\\
\indent Now let $p(T)(t) = \sum_{r \geq 0} p_{r}(T) t^{r}$ denote the
Pontryagin polynomial. The relationship between the Pontryagin classes
and the Chern classes is given by
\[
p(T)(-1) = c(T)(1) \cdot c(T)(-1) \indent \cite[\mbox{Cor. } 15.5]{MS}.
\]
Hence $p(T)(-1)$ equals
\[
\prod_{k \neq l} (1+ \delta_{k} -\delta_{l})^{2\g} =
\prod_{k<l} (1 - (\delta_{k} - \delta_{l})^{2})^{2\g}.
\]
The total Pontryagin class of $\mnd$ then equals $p(T)(1)$ and hence
the result. $\Box$
\section{A Complete Set of Relations.}
Whilst we observed in remark \ref{inadequacy} that neither the Mumford relations nor the
dual
Mumford relations are in themselves a complete set of relations when
the rank is greater than two, it is still
possible to put these relations into the context of the 
completeness criteria. In terms of
these criteria we will show how
the Mumford relations contain subsets corresponding to all strata of
the form
\[
\mu= (d_{1}/n_{1},...,d_{P}/n_{P})
\]
where $n_{P}=1.$ Similarly the dual Mumford relations contain subsets
corresponding to all
those strata with $n_{1}=1.$ From this we shall deduce that in the
rank three case the Mumford and dual Mumford relations form a complete
set.\\
\indent Before we continue with the main proposition we need a lemma on the
vanishing of the Mumford and dual Mumford
relations on restriction to a stratum.
\begin{lem}
\label{vanishing}
 Let $\mu=(d_{1}/n_{1},...,d_{P}/n_{P})$. The image of
the Mumford relation $\sigma_{r,S}^{k}$ under the restriction map
\[
H^{*}_{\cal G}({\cal C};{\bf Q}) \rightarrow H^{*}_{\cal G}(\Cmu;{\bf
Q})
\]
vanishes when $r < d_{P}/n_{P} -2g +1$. The image of the dual Mumford
relation $\tau_{r,S}^{k}$ under the restriction map vanishes when $r <
2\g-d_{1}/n_{1}.$
\end{lem}
{\bf Proof}
Recall that the Mumford relations are given by
$\sigma_{r,S}^{k} (r<0,0 \leq k \leq n-1,S \subseteq \{1,...,2g\})$
when $\Psi(t)= t^{d-n\g}c(\pi_{!}{\cal V})(t^{-1})$ is written in the form
\[
\sum_{r=-\infty}^{\g} (\sigma_{r}^{0} + \sigma_{r}^{1} t + \cdot \cdot
\cdot + \sigma_{r}^{n-1} t^{n-1})(\om(t))^{r}, \indent
\sigma_{r}^{k}=\sum_{S \subseteq \{1,...,2g\}} \sigma_{r,S}^{k}
\prod_{s \in S} b_{1}^{s}.
\]
For $1 \leq k \leq n$ and any fixed integer $R$ the power $t^{-k}$ appears in
\[
\sum_{r=-\infty}^{\g} (\sigma_{r}^{0} + \sigma_{r}^{1} t + \cdot \cdot \cdot +
\sigma_{r}^{n-1} t^{n-1})(\om(t))^{r-R-1}
\]
only when $r=R.$ Let $C_{r}^{i}$ denote the
coefficient of $t^{-i}$ in $\Psi(t)(\tilde{\Omega}(t))^{-r-1}$. Then
\[
(\sigma_{r}^{0} + \sigma_{r}^{1} t + \cdot \cdot \cdot +
\sigma_{r}^{n-1} t^{n-1}) = (t^{n} + a_{1} t^{n-1} + \cdot \cdot \cdot
+ a_{n})  \sum_{i=1}^{n} C_{r}^{i} t^{-i}
\]
modulo negative powers of $t$ and hence
\begin{equation}
\sigma_{r}^{n-k} = \sum_{i=1}^{k} a_{k-i} C_{r}^{i} \indent (r<0,1
\leq k \leq n) \label{31}.
\end{equation}
\indent Now let $K$ be a fixed line bundle over $M$ of degree $D$ where $D$ is
the smallest integer such that
\[
\mu(Q_{P} \otimes K) = \frac{d_{P}}{n_{P}} + D > 2\g
\]
where $Q_{P}=E_{P}/E_{P-1}.$ Since $\mu(Q_{p} \otimes K) \geq
\mu(Q_{P} \otimes K) > 2\g$ then $\pi_{!}({\cal V}_{p} \otimes
\phi^{*}K)$ is a bundle over ${\cal C}(n_{p},d_{p})^{ss}$ of rank
$d_{p} +(D-\g)n_{p}$ for each $1 \leq p \leq P.$ In particular
\[
\Psi(\pi_{!}({\cal V}_{p} \otimes \phi^{*}K))(t) =
t^{d_{p}+n_{p}(D-\g)} c (\pi_{!}({\cal V}_{p} \otimes \phi^{*}K))(t^{-1})
\]
is a polynomial modulo relations in $H^{*}_{{\cal G}(n_{p},d_{p})}
({\cal C}(n_{p},d_{p})^{ss};{\bf Q})$. From GRR we have that
$\ch(\pi_{!}({\cal V}_{p} \otimes \phi^{*}K))$ equals
\begin{equation}
 \ch(\pi_{!}{\cal V}_{p}) + \pi_{*}(\ch {\cal V}_{p} \cdot 1 \otimes
D\omega) = \ch(\pi_{!}{\cal V}_{p}) +D  \sum_{k=1}^{n_{p}}
e^{\delta_{k}^{p}}. \label{32}
\end{equation}
In terms of Chern polynomials (\ref{32}) gives
\[
c(\pi_{!}({\cal V}_{p} \otimes \phi^{*}K))(t) =
(\Omega_{p}(t))^{D} c(\pi_{!}{\cal V}_{p})(t)
\]
where $\Omega_{p}(t) = \prod_{k=1}^{n_{p}} ( 1 + \delta_{k}^{p}t)$. Hence
\begin{equation}
\prod_{p=1}^{P} \Psi(\pi_{!}({\cal V}_{p} \otimes \phi^{*}K))(t)
= (\tilde{\Omega}_{\mu}(t))^{D} \Psi_{\mu}(t) \label{51}
\end{equation}
is a polynomial modulo relations in $\HG(\Cmu;\Q)$ where
$\Psi_{\mu}(t),$ and $\om_{\mu}(t)$ are respectively the restrictions
to $\HG(\Cmu;\Q)$ of $\Psi(t)$ and $\om(t)$. Thus the coefficient of $
t^{-k}$ in
$\Psi_{\mu}(t)\tilde{\Omega}_{\mu}(t)^{-r-1}$ is a relation
when $r \leq -1 -D$. So by (\ref{31}) the restriction of $\sigma_{r}^{k}$ to
$H^{*}_{\cal
G}(\Cmu;{\bf Q})$ vanishes when $r \leq d_{P}/n_{P} - 2g.$ The dual
calculation follows by a similar argument. $\indent \Box$\\[\baselineskip]
\indent Thus finally we come to
\begin{prop}
\label{biggy}
Let $\mu = (d_{1}/n_{1},...,d_{P}/n_{P})$ with
$n_{P}=1$. Then there is a subset ${\cal R}_{\mu}$ of the ideal
generated by the Mumford relations such that the image of the ideal
generated by ${\cal R}_{\mu}$ under the restriction map
\[
\HG(\C;\Q) \rightarrow \HG(\Cnu;\Q) \indent
\nu=(\tilde{d}_{1}/\tilde{n}_{1},...,\tilde{d}_{T}/\tilde{n}_{T})
\]
is zero when either
\[
\mbox{(i) } \tilde{d}_{T}/\tilde{n}_{T}> d_{P} \indent \mbox{or}
\indent \mbox{(ii) } \tilde{n}_{T}=1, \tilde{d}_{T} = d_{P}, \mbox{
and } \nu \not \geq \mu
\]
and contains the ideal of $\HG(\Cmu;\Q)$ generated by $e_{\mu}$ when
$\nu=\mu.$\\
\indent Let $\mu = (d_{1}/n_{1},...,d_{P}/n_{P})$ with
$n_{1}=1$. Then there is a subset ${\cal R}_{\mu}$ of the ideal
generated by the dual Mumford relations such that the image of the ideal
generated by ${\cal R}_{\mu}$ under the restriction map
\[
\HG(\C;\Q) \rightarrow \HG(\Cnu;\Q) \indent
\nu=(\tilde{d}_{1}/\tilde{n}_{1},...,\tilde{d}_{T}/\tilde{n}_{T})
\]
is zero when either
\[
\mbox{(i) } \tilde{d}_{1}/\tilde{n}_{1}< d_{1}/n_{1} \indent \mbox{or}
\indent \mbox{(ii) } \tilde{n}_{1}=1, \tilde{d}_{1} = d_{1} \mbox{
and } \nu \not \geq \mu
\]
and contains the ideal of $\HG(\Cmu;\Q)$ generated by $e_{\mu}$ when
$\nu=\mu.$
\end{prop}
{\bf Proof}
Let $\Psi(t)= t^{d-n\g}c(\pi_{!}{\cal V})(t^{-1})$ and
let $C^{R}_{K}, (R<0,1 \leq K \leq n)$ denote the coefficient of
$t^{-K}$ in $\Psi(t)(\om(t))^{-R-1}$. Let
\[
\mu=(d_{1}/n_{1},...,d_{P-1}/n_{P-1},d_{P})
\]
so that $n_{P}=1.$\\
\indent Since the Chern polynomial is multiplicative the restriction
in $\HG(\Cmu;\Q)$ of $C_{R}^{K}$, which we will write
as $C_{R}^{K,\mu}$, equals the coefficient of $t^{-1}$ in
\begin{equation}
t^{K-1} \prod_{p=1}^{P} \Psi_{p}(t)(\om_{p}(t))^{-R-1} \label{50}
\end{equation}
where
\[
\Psi_{p}(t)= t^{d_{p}-n_{p}\g}c(\pi_{!}{\cal V}_{p})(t^{-1}), \indent
\om_{p}(t)=t^{n_{p}}+a_{1}^{p} t^{n_{p}-1}+ \cdot \cdot \cdot +
a_{n_{p}}^{p}
\]
for $1 \leq p \leq P$. Further from the previous lemma we know that
$C_{R}^{K,\mu}$ vanishes when $R<-D=d_{P}-2g+1.$\\
\indent We facilitate the proof of proposition \ref{biggy} with the following
lemma and corollaries
\begin{lem}
Let $\theta(t)$ equal
\begin{equation}
t^{d-nd_{P}+(n-1)\g} \prod_{p=1}^{P-1} \prod_{k=1}^{n_{p}} ( 1
+(\delta_{k}^{p}-a_{1}^{P})/t)^{W_{k}^{p}+\g-d_{P}} \exp \left\{
\frac{\Xi_{k,1}^{p,P}}{t+\delta_{k}^{p}-a_{1}^{P}} \right\}. \label{52}
\end{equation}
Then modulo relations in $\HG(\Cmu;\Q)$,
\[
C_{-D}^{K,\mu}= (-a_{1}^{P})^{K-1} (\xi_{1,1}^{P,P})^{g} \Theta
\]
where $\Theta$ is the constant coefficient of $\theta(t).$
\end{lem}
{\bf Proof}
From corollary \ref{Deltaresn} we know that
\[
\Psi_{P}(t) (\om_{P}(t))^{D-1} = (t+a_{1}^{P})^{\g} \exp \left\{
\frac{ \xi_{1,1}^{P,P}}{t+a_{1}^{P}} \right\}
\]
where $\xi_{1,1}^{P,P} = \sum_{s=1}^{g} b_{1}^{P,s} b_{1}^{P,s+g}.$
Also in a Laurent series the coefficient of $t^{-1}$ is invariant under
transformations such as $t \mapsto t-a_{1}^{P}.$ So from (\ref{50})
$C_{-D}^{K,\mu}$ equals the coefficient of $t^{-1}$ in
\begin{equation}
(t-a_{1}^{P})^{K-1} t^{\g} \exp (\xi_{1,1}^{P,P}/t) \prod_{p=1}^{P-1}
\Psi_{p}(t-a_{1}^{P})(\om_{p}(t-a_{1}^{P}))^{D-1}. \label{D}
\end{equation}
\indent From the proof of lemma \ref{vanishing} (\ref{51}) we know that
\[
\Psi_{p}(t)(\om_{p}(t))^{D-1} = \Psi(\pi_{!}({\cal V}_{p} \otimes
\phi^{*}{\cal L}))(t)
\]
where ${\cal L}$ is a fixed line bundle over $M$ of degree $D-1.$ For each
$p \neq P$, $Q_{p} \otimes {\cal L}$ is a semistable bundle of
slope
\[
\frac{d_{p}}{n_{p}} - d_{P} + 2\g > 2\g.
\]
Hence $\pi_{!}({\cal
V}_{p} \otimes \phi^{*}{\cal L})$ is a bundle over
$\C(n_{p},d_{p})^{ss}$ and $\Psi_{p}(t)(\om_{p}(t))^{D-1}$ is a polynomial
modulo relations in
$\HS_{\G (n_{p},d_{p})}(\C(n_{p},d_{p})^{ss};\Q).$ As
$(\xi_{1,1}^{P.P})^{g+1} = 0$ it follows from (\ref{D}) that
$C_{-D}^{K,\mu}$ equals the constant coefficient of
\begin{equation}
(\xi_{1,1}^{P,P})^{g} (t-a_{1}^{P})^{K-1} \prod_{p=1}^{P-1}
\Psi_{p}(t-a_{1}^{P}) (\om_{p}(t-a_{1}^{P}))^{D-1} \label{E}
\end{equation}
modulo relations in $\HG(\Cmu;\Q)$.\\
\indent Since $\sum_{k=1}^{n_{p}}
W_{k}^{p} =d_{p}$ then we know from corollary \ref{Chernresn} that $\Psi_{p}(t-a_{1}^{P})$
equals
\[
(\om_{p}(t-a_{1}^{P}))^{-\g} t^{d_{p}} \prod_{k=1}^{n_{p}} ( 1 +(
\delta_{k}^{p} -a_{1}^{P})/t)^{W_{k}^{p}} \exp \left\{
\frac{X_{k}^{p}}{t+\delta_{k}^{p} -a_{1}^{P}} \right\}.
\]
Recall from lemma \ref{normal} that
\[
\Xi_{k,1}^{p,P} = \sum_{s=1}^{g} \left( \sum_{i=1}^{n_{p}} b_{i}^{p,s}
\frac{\partial \delta_{k}^{p}}{\partial a_{i}^{p}} - b_{1}^{P,s}
\right) \left( \sum_{i=1}^{n_{p}} b_{i}^{p,s+g}
\frac{\partial \delta_{k}^{p}}{\partial a_{i}^{p}} - b_{1}^{P,s+g}
\right)
\]
and we also have that
\[
X_{k}^{p} =  \sum_{s=1}^{g} \left( \sum_{i=1}^{n_{p}} b_{i}^{p,s}
\frac{\partial \delta_{k}^{p}}{\partial a_{i}^{p}}
\right) \left( \sum_{i=1}^{n_{p}} b_{i}^{p,s+g}
\frac{\partial \delta_{k}^{p}}{\partial a_{i}^{p}}
\right).
\]
Since
\[
(\xi_{1,1}^{P,P})^{g} = (-1)^{g\g/2} g! \prod_{s=1}^{2g}
b_{1}^{P,s}
\]
then
\[
(\xi_{1,1}^{P,P})^{g}(\Xi_{k,1}^{p,P})^{q} =
(\xi_{1,1}^{P,P})^{g}(X_{k}^{p})^{q} \indent (q \geq 0).
\]
Thus by (\ref{E}) and the identity $\om_{p}(t-a_{1}^{P}) =
t^{n_{p}}\prod_{k=1}^{n_{p}}(1+(\delta_{k}^{p}-a_{1}^{P})/t),$ we have
that $C_{-D}^{K,\mu}$ equals the constant coefficient of
\[
(\xi_{1,1}^{P,P})^{g} (t-a_{1}^{P})^{K-1} \theta(t).
\]
Since $(\xi_{1,1}^{P,P})^{g}\theta(t)$ is a polynomial modulo
relations in $\HG(\Cmu;\Q)$ then the lemma follows. $\indent
\Box.$
\begin{cor}
Define $C_{R,S}^{K} (R<0,1 \leq K \leq n, S \subseteq
\{1,...,2g\})$ by
\[
C_{R}^{K} = \sum_{S \subseteq \{1,...,2g\}} C_{R,S}^{K} \prod_{s \in
S} b_{1}^{s}
\]
writing $C_{R,S}^{K}$ in terms of the elements (\ref{1000})
and also define $\ta_{r},\tb_{r}^{s}$ and $\tf_{r}$ by
\[
c_{r}(\bigoplus_{p=1}^{P-1} {\cal V}_{p}) = \ta_{r} \otimes 1 +
\sum_{s=1}^{2g} \tb_{r}^{s} \otimes \alpha_{s} + \tf_{r} \otimes \omega.
\]
Then the restriction of $C_{-D,S}^{K}$ to $\HG(\Cmu;\Q)$ equals a
non-zero constant multiple of
\begin{equation}
(a_{1}^{P})^{K-1} \prod_{s \not \in S} (\tb_{1}^{s}-(n-1) b_{1}^{P,s})
\Theta \label{M}
\end{equation}
for any subset $S \subseteq \{1,...,2g\}.$
\end{cor}
{\bf Proof}
We know that $(\xi_{1,1}^{P,P})^{g}$ equals
\[
(-1)^{g\g/2} g! \prod_{s=1}^{2g} b_{1}^{P,s} = (-1)^{g\g/2}
n^{-2g} g! \prod_{s=1}^{2g} ((\tb_{1}^{s} +
b_{1}^{P,s})-(\tb_{1}^{s}-(n-1) b_{1}^{P,s}))
\]
and also that the restriction of $b_{1}^{s}$ in $\HG(\Cmu;\Q)$ equals
$\tb_{1}^{s}+ b_{1}^{P,s}.$ Further
\[
\tb_{1}^{s} - (n-1)b_{1}^{P,s} = \sum_{p=1}^{P-1} \sum_{k=1}^{n_{p}} \left(
\sum_{i=1}^{n_{p}} b_{1}^{p,s} \frac{\partial \delta_{k}^{p}}{\partial
a_{i}^{p}}  - b_{1}^{P,s} \right).
\]
So the corollary follows once we note from
(\ref{52}) that $\theta(t)$, and hence $\Theta$, can be written in
terms of the elements (\ref{NEW1}). $\indent \Box$.
\begin{cor}
Let $\Lambda$ equal
\begin{equation}
\bigcup \{ \sigma_{-D,S}^{n-1},...,\sigma_{-D,S}^{0} \} \label{54}
\end{equation}
where the union varies over all subsets $S \subseteq \{1,...,2g\}$.
Then all elements of the form
\begin{equation}
\prod_{k=2}^{n-1} (\tf_{k})^{m_{k}} \prod_{k=1}^{n-1} \prod_{s \in
S_{k}} \tb_{k}^{s} \prod_{k=1}^{n-1} (\ta_{k})^{r_{k}} (a_{1}^{P})^{r}
\prod_{s \in S} b_{1}^{P,s} \Theta \label{81}
\end{equation}
lie in the restriction of the ideal generated by $\Lambda$, where
$r,r_{1},...,r_{n-1},m_{2},...,m_{n-1}$ are arbitrary non-negative
integers and $S,S_{1},...,S_{n-1}$ are subsets of
$\{1,...,2g\}$.
\end{cor}
{\bf Proof}
Let $(\Lambda)$ denote the ideal of $\HG(\C;\Q)$
generated by $\Lambda$. Using induction on (\ref{31}) we know that the
restriction of $C_{-D,S}^{K}$ lies in the image of
$(\Lambda)$. From (\ref{M}) and since $b_{1}^{s}$ restricts to $\tb_{1}^{s}
+ b_{1}^{P,s}$ it follows that all elements of the form
\[
(a_{1}^{P})^{K-1} \prod_{s \in S_{1}} \tb_{1}^{s} \prod_{s \in S_{2}}
b_{1}^{P,s} \Theta
\]
for arbitrary $S_{1},S_{2} \subseteq \{1,...,2g\}$ and $1 \leq K \leq
n$, lie in the restriction of $(\Lambda).$ The restriction of
$a_{k}$ in $\HG(\Cmu;\Q)$ equals
$\ta_{k}+\ta_{k-1}a_{1}^{P}$. By noting that $(a_{1}^{P})^{r}$ equals
\[
(\ta_{1}+a_{1}^{P})(a_{1}^{P})^{r-1} -(\ta_{2} +
\ta_{1} a_{1}^{P})(a_{1}^{P})^{r-2} + \cdot \cdot \cdot +
(-1)^{n-1}(\ta_{n-1}a_{1}^{P})(a_{1}^{P})^{r-n}
\]
for $r \geq n$ we see that all elements of the form
\[
(a_{1}^{P})^{r} \prod_{s \in S_{1}} \tb_{1}^{s} \prod_{s \in S_{2}}
b_{1}^{P,s} \cdot \Theta \indent (r \geq 0)
\]
lie in the restriction of $(\Lambda)$. Finally working inductively on
the variables $r_{1},...,r_{n-1},$ $S_{2},S_{3},...,S_{n-1}$ and
$m_{2},m_{3},...m_{n-1}$ in that order we find that all elements of
the form (\ref{81}) lie in the image of $(\Lambda)$ since under the
restriction map $\HG(\C;\Q) \rightarrow \HG(\Cmu;\Q)$
\[
a_{k} \mapsto \ta_{k}+\ta_{k-1}a_{1}^{P} \indent b_{k}^{s} \mapsto
\tb_{k}^{s}+ a_{1}^{P} \tb_{k-1}^{s} + \ta_{k-1}b_{1}^{P,s}
\]
and
\begin{equation}
f_{k} \mapsto \tf_{k} + d_{P}\ta_{k-1} + a_{1}^{P} \tf_{k-1} +
\sum_{s=1}^{g} ( \tb_{k-1}^{s} b_{1}^{P,s+g} + b_{1}^{P,s}
\tb_{k-1}^{s+g}). \quad \Box \label{60}
\end{equation}
\indent We now continue with the proof of proposition \ref{biggy}. Let
$\C'=\C(n-1,d-d_{P})$ and let $\G'= \G(n-1,d-d_{P}).$
Let $\mu' = (d_{1}/n_{1},...,d_{P-1}/n_{P-1})$ and let $e_{\mu'}$ denote the
equivariant Euler class of the normal bundle to $\C'_{\mu'}$ in $\C'.$ Let
\[
U_{\mu'} = \C' - \bigcup_{\nu'>\mu'} \C'_{\nu'}.
\]
Then $U_{\mu'}$ is an open subset of $\C'$ which contains
$\C'_{\mu'}$ as a closed submanifold. So we have the maps\\
\begin{picture}(400,120)
\put(190,60){\makebox(0,0){$H^{*-2d_{\mu'}}_{\G'}(\C'_{\mu'};\Q) \rightarrow
H^{*}_{\G'}(U_{\mu'};\Q) \rightarrow
H^{*}_{\G'}(U_{\mu'}-\C'_{\mu'};\Q)$}}
\put(190,110){\makebox(0,0){$H^{*}_{\G'}(\C';\Q)$}}
\put(190,10){\makebox(0,0){$H^{*}_{\G'}(\C'_{\mu'};\Q)$}}
\put(190,102){\vector(0,-1){34}}
\put(190,50){\vector(0,-1){32}}
\put(100,50){\vector(2,-1){65}}
\put(150,30){\makebox(0,0)[tr]{multiplication by $e_{\mu'}$}}
\end{picture}
\\
Let $a'_{r},{b^{s}_{r}}'$ and $f'_{r}$ denote the generators of
$H^{*}_{\G'}(\C';\Q)$. Also take $\nu' \not \geq \mu'$ and let
$\tta_{r},\ttb_{r}^{s},\ttf_{r}$ denote the restrictions of
$a'_{r},{b^{s}_{r}}',f'_{r}$ in $H^{*}_{\G'}(\C'_{\nu'};\Q)$. Since the
stratification is equivariantly perfect then the restriction map
\[
H^{*}_{\G'}(\C';\Q) \rightarrow H^{*}_{\G'}(U_{\mu'};\Q)
\]
is surjective \cite[p.859]{K2}. From the exactness of the Thom-Gysin
sequence we have that for every element of
the form $\alpha e_{\mu'}$ in $H^{*}_{\G'}(\C'_{\mu'};\Q)e_{\mu'}$ there is
some $\beta
(a'_{r},{b_{r}^{s}}' ,f'_{r})$ in $\HS_{\G'}(\C';\Q)$ such that
\[
\beta(\ta_{r},\tb_{r}^{s},\tf_{r}) = \alpha e_{\mu'} \mbox{ and }
\beta(\tta_{r},\ttb_{r}^{s},\ttf_{r}) = 0.
\]
Since every element of the form (\ref{81}) lies in the restriction of
$(\Lambda)$ to $\HG(\Cmu;\Q)$ then every element of the form
\begin{equation}
\beta(\ta_{r},\tb_{r}^{s},\tf_{r}) (a_{1}^{P})^{r} \prod_{s \in S}
b_{1}^{P,s} \Theta \indent (r \geq 0,S \subseteq \{1,...,2g\}) \label{z1}
\end{equation}
similarly lies in the restriction of $(\Lambda)$. Now let
$\nu=(\nu',d_{P})$ with $\nu' \not \geq \mu'$. Note that the
restriction map
\[
\HG(\C;\Q) \rightarrow \HG(\Cnu;\Q)
\]
is formally the same as (\ref{60}) but with
$\tta_{r},\ttb_{r}^{s},\ttf_{r}$ replacing
$\ta_{r},\tb_{r}^{s},\tf_{r}$. Thus there are elements of $(\Lambda)$
which restrict to (\ref{z1}) under (\ref{60}) and have restriction
\[
\beta(\tta_{r},\ttb_{r}^{s},\ttf_{r}) (a_{1}^{P})^{r} \prod_{s \in S}
b_{1}^{P,s} \hat{\Theta} = 0
\]
in $\HG(\Cnu;\Q)$.\\[\baselineskip]
\indent Define ${\cal R}_{\mu}$ to be all those elements
of $(\Lambda)$ which restrict to an element of the form
\[
\alpha e_{\mu'} (a_{1}^{P})^{r} \prod_{s \in S} b_{1}^{P,s} \Theta
\indent (r \geq 0, S \subseteq \{1,...,2g\}, \alpha \in
\HS_{\G'}(\C'_{\mu'};\Q))
\]
in $\HG(\Cmu;\Q)$ and which restrict to zero in $\HG(\Cnu;\Q)$ for any
$\nu = (\nu',d_{P})$ with $\nu' \not \geq \mu'$.\\[\baselineskip]
\indent From the definition of $\Theta$ (\ref{52}) we know that
$e_{\mu'} \Theta$ is the constant coefficient of
\begin{equation}
(-1)^{d_{\mu'}} t^{d_{\mu'}} c({\cal N}_{\mu'})(-t^{-1}) \theta(t) \label{53}
\end{equation}
where ${\cal N}_{\mu'}$ is the normal bundle to $\C'_{\mu'}$ in
$\C'$ and $d_{\mu'}$ is the codimension of $\C'_{\mu'}$ in
$\C'$. From lemma \ref{normal} and the fact that
\[
d_{\mu'} + d - nd_{P} + (n-1)\g = d_{\mu}
\]
we know (\ref{53}) equals
\[
(-1)^{d_{\mu'}} t^{d_{\mu}} c({\cal N}_{\mu})(-t^{-1})
\]
which has constant coefficient $(-1)^{d_{\mu'}+d_{\mu}}e_{\mu}.$ Hence the
ideal
\[
\HG(\Cmu;\Q)e_{\mu}
\]
lies in the restriction of ${\cal R}_{\mu}$ to
$\HG(\Cmu;\Q).$\\
\indent Finally from lemma \ref{vanishing} and the definition of $\Lambda$ (\ref{54})
we know that the image of ${\cal R}_{\mu}$ under the restriction map
\[
\HG(\C;\Q) \rightarrow \HG(\Cnu;\Q) \indent \nu =
(\tilde{d}_{1}/\tilde{n}_{1},...,\tilde{d}_{T}/\tilde{n}_{T})
\]
vanishes when $\tilde{d}_{T}/\tilde{n}_{T} > d_{P}/n_{P}$ proving the
first half of proposition \ref{biggy}.\\
\indent The proof of the dual case follows in a similar fashion.$
\indent \Box$\\[\baselineskip]
\indent In the general rank case there are strata of types not covered in the
previous proposition. Moreover the strata on which the restrictions of
the relations have been demonstrated to vanish do not generally coincide with
the strata mentioned
in the hypotheses of the completeness criteria. However in the
rank two and rank three cases all unstable strata are covered by the above
proposition. In the rank two case proposition \ref{biggy} shows that the Mumford
relations and the dual Mumford relations both form complete sets,
simply duplicating Kirwan's work \cite{K2} and remark \ref{dualise}. In the
rank three case we have the following:\\[\baselineskip]
{\bf THEOREM 1.} {\em The Mumford and dual Mumford relations together with the
relation (\ref{NORM}) due to the normalisation of the
universal bundle $V$ form a complete set of relations for $H^{*}({\cal
M}(3,d);\Q).$}
{\bf Proof}
The unstable strata are now of types (2,1),(1,1,1) and
(1,2). From the previous proposition we may meet the completeness
criteria for the (2,1) and (1,1,1) strata using the Mumford relations.
In these cases those strata where the
restriction of ${\cal R}_{\mu}$ have been shown to vanish are those strata
$\Cnu$ such that $\nu \prec \mu$. The criteria for the (1,2) types may be met
using
the dual Mumford relations. In this case those strata where the
restriction of ${\cal R}_{\mu}$ vanishes (according to proposition \ref{biggy}) are those strata $\Cnu$ such that $\nu \not \geq \mu$ which certainly
includes those strata such that $\nu \prec \mu.
\indent \Box$
\begin{rem}
As remarked earlier it was shown in \cite[thm.4]{E} that the Mumford
relations $\sigma_{-1,S}^{1}$ for $S \subseteq \{1,...,2g\}$ generate
the relation ideal of $\HS({\cal M}_{0}(2,1);\Q)$ as a
$\Q[a_{2},f_{2}]$-module. Evidence for this theorem appears in the
Poincar\'{e} polynomial of the relation ideal which equals \cite[p.593]{AB}
\[
\frac{t^{2g}(1+t)^{2g}}{(1-t^{2})(1-t^{4})}.
\]
\indent Similarly in the rank three case the Poincar\'{e} polynomial
of the ideal of relations among our generators for $H^{*}({\cal
M}_{0}(3,1);\Q)$ equals
\[
\frac{(1+t^{2})^{2} t^{4g-2}(1+t)^{2g}(1+t^{3})^{2g} - (1+t^{2}+t^{4})
t^{6g-2}(1+t)^{4g}}{(1-t^{2})(1-t^{4})^{2}(1-t^{6})},
\]
The first Mumford relation $\sigma_{-1,\{1,...,2g\}}^{2}$ has degree $4g-2$
and the first dual Mumford relation $\tau_{-1,\{1,...,2g\}]}^{2}$ has degree
$4g$. This strongly suggests that the relations
\[
\{ \sigma_{-1,S}^{i}, \tau_{-1,S}^{i} : i=1,2, S \subseteq \{1,...,2g\} \}
\]
generate the relation ideal of $\HS({\cal N}(3,d);\Q)$ as a
\[
\Q[a_{2},a_{3},f_{2},f_{3}] \otimes
\Lambda^{*}\{b_{2}^{1},...,b_{2}^{2g}\}
\]
module.
\end{rem}
\section{On the Vanishing of the Pontryagin Ring.}
\indent We now move on to discuss the Pontryagin ring of the moduli
space in the rank three case. For each $S \subseteq \{1,...,2g\}$ we
define $\Psi_{S}(t)$ and $\Psi^{*}_{S}(t)$ by writing
\[
\Psi(t) = \sum_{S \subseteq \{1,...,2g\}} \Psi_{S}(t) \prod_{s \in S}
b_{1}^{s}, \indent \Psi^{*}(t) = \sum_{S \subseteq \{1,...,2g\}}
\Psi^{*}_{S}(t) \prod_{s \in S} b_{1}^{s}.
\]
Kirwan proved the Newstead-Ramanan conjecture \cite[$\S$ 4]{K2} by considering
relations derived from the expression
\[
\Psi_{\{1,...,2g\}}(t)\Psi_{\{1,...,2g\}}(-t-a_{1}).\]
Arguing along similar lines but now considering the expression
\[
\Phi(t) = \Psi_{\{1,...,2g\}}(t)\Psi^{*}_{\{1,...,2g\}}(t)
\]
we will show that in the rank three case the Pontryagin ring vanishes
in degree $12g-8$ and above -- theorem 2 below.
\begin{lem}
\label{Pontlemma}
Let $\mu = (d_{1}, d_{2},...,d_{n}) \in \Delta$. The
restriction of $\Phi(t)$ to $\HG(\Cmu;\Q)$ equals
\[
(-1)^{g} \frac{A(t)^{2g}}{n^{4g}\om_{\mu}(t)}
\]
where
\[
\om_{\mu}(t) = \prod_{p=1}^{n} (t+a_{1}^{p}), \indent A(t) =
\sum_{p=1}^{n} \prod_{q \neq p} (t+a_{1}^{q}).
\]
\end{lem}
{\bf Proof}
From corollary \ref{Deltaresn} we know that the restriction of
$\Psi(t)$ to $\HG(\Cmu;\Q)$ equals
\[
\prod_{p=1}^{n} (t+a_{1}^{p})^{d_{p}-\g} \exp \left\{
\frac{\xi_{p}}{t+a_{1}^{p}} \right\}
\]
where $\xi_{p} = \xi_{1,1}^{p,p} = \sum_{s=1}^{g} b_{1}^{p,s} b_{1}^{p,s+g}.$
Let $v_{s} = b_{1}^{1,s} + \cdot
\cdot \cdot + b_{1}^{n,s}$ denote the restriction of $b_{1}^{s}$ to
$\HG(\Cmu;\Q)$ and let $w_{i,j}^{s} = b_{1}^{i,s} - b_{1}^{j,s}$ (see
(\ref{1001})). Then
$nb_{1}^{i,s} = v_{s} + \sum_{j=1}^{n} w_{i,j}^{s}$ and hence
\[
n^{2} \xi_{i} = \sum_{s=1}^{g} v_{s}v_{s+g} + \sum_{s=1}^{g} \left(
v_{s} \sum_{j=1}^{n} w_{i,j}^{s+g} + \sum_{j=1}^{n} w_{i,j}^{s}
v_{s+g} \right) + \sum_{s=1}^{g} \sum_{j=1}^{n} \sum_{k=1}^{n} w_{i,j}^{s}
w_{i,k}^{s+g}.
\]
Note that
\begin{equation}
\sum_{p=1}^{n} \frac{\xi_{p}}{t+a_{1}^{p}} = \frac{1}{\om_{\mu}(t)}
\sum_{i=1}^{n} \sum_{q \neq i} \xi_{i} (t+a_{1}^{q}) \label{N}
\end{equation} Thus (\ref{N}) equals
\[
\frac{1}{n^{2} \om_{\mu}(t)} \left\{ A(t) \sum_{s=1}^{g} v_{s}v_{s+g}
+ \sum_{s=1}^{g} (B_{s}(t) v_{s+g} + v_{s}B_{s+g}(t)) + \Gamma(t)
\right\}
\]
where
\[
A(t) = \sum_{i=1}^{n} \prod_{q \neq i} (t+a_{1}^{q}) , \indent B_{s}(t)
= \sum_{i=1}^{n} \sum_{j=1}^{n} w_{i,j}^{s} \prod_{q \neq i} 
(t+a_{1}^{q}),
\]
\[
\Gamma(t) = \sum_{i=1}^{n} \sum_{j=1}^{n}
\sum_{k=1}^{n} \sum_{s=1}^{g} w_{i,j}^{s} w_{i,k}^{s+g} \prod_{q \neq i}(t+a_{1}^{q}).
\]
The exponential of (\ref{N}) equals
\[
\exp \left\{ \frac{\Gamma(t)}{n^{2}\om_{\mu}(t)} \right\}
\prod_{s=1}^{g} \left[ 1+ \frac{B_{s}(t)v_{s+g} + v_{s}
B_{s+g}}{n^{2}\om_{\mu}(t)} + \left( \frac{A(t)}{n^{2}\om_{\mu}(t)} -
\frac{B_{s}B_{s+g}}{n^{4}\om_{\mu}(t)^{2}}\right) v_{s}v_{s+g} \right].
\]
The coefficient of $\prod_{s=1}^{2g} v_{s}$ in the above then equals
\[
(-1)^{g\g/2} \exp \left\{ \frac{\Gamma(t)}{n^{2}\om_{\mu}(t)} \right\}
\prod_{s=1}^{g} \left( \frac{A(t)}{n^{2}\om_{\mu}(t)} -
\frac{B_{s}B_{s+g}}{n^{4}\om_{\mu}(t)^{2}} \right)
\]
or equivalently
\[
(-1)^{g\g/2} \exp \left\{ \frac{\Gamma(t)}{n^{2}\om_{\mu}(t)} \right\}
\left( \frac{A(t)}{n^{2}\om_{\mu}(t)} \right)^{g}
\exp\left\{\frac{-\xi(t)}{n^{2}A(t)\om_{\mu}(t)} \right\}
\]
where $\xi(t) = \sum_{s=1}^{g} B_{s}(t)B_{s+g}(t)$. Thus the
restriction of $\Psi_{\{1,...,2g\}}(t)$ to $\HG(\Cmu;\Q)$ equals
\[
(-1)^{g\g/2}\left(\prod_{p=1}^{n} (t+a_{1}^{p})^{d_{p}-\g} \right)
\exp \left\{ \frac{\Gamma(t)}{n^{2}\om_{\mu}(t)} \right\}
\left( \frac{A(t)}{n^{2}\om_{\mu}(t)} \right)^{g}
\exp\left\{\frac{-\xi(t)}{n^{2}A(t)\om_{\mu}(t)} \right\}
\]
and similarly the restriction of $\Psi^{*}_{\{1,...,2g\}}(t)$ to
$\HG(\Cmu;\Q)$ equals
\[
(-1)^{g\g/2} \left(\prod_{p=1}^{n} (t+a_{1}^{p})^{3\g+1-d_{p}} \right)
\exp \left\{ \frac{- \Gamma(t)}{n^{2}\om_{\mu}(t)} \right\}
\left(\frac{-A(t)}{n^{2}\om_{\mu}(t)} \right)^{g}
\exp\left\{\frac{\xi(t)}{n^{2}A(t)\om_{\mu}(t)} \right\}.
\]
The result then follows. $\indent \Box$\\[\baselineskip]
\indent Now if we write $\Phi(t)$ in the form
\[
\sum_{r = -\infty}^{2g-1} ( \rho_{r}^{0} + \rho_{r}^{1} t + \cdot
\cdot \cdot + \rho_{r}^{n-1} t^{n-1}) ( \om(t))^{r}
\]
where $\om(t) = t^{n} + a_{1}t^{n-1} + \cdot \cdot \cdot +a_{n}$ then
we know that the elements $\rho_{r}^{k},(r < 0 ,0 \leq k \leq n-1)$
lie in the kernel of the restriction map
\[
\HG(\C;\Q) \rightarrow \HG(\Css;\Q).
\]
\indent From lemma \ref{Pontlemma} we know that the restriction of $\Phi(t)$ to
$\HG(\Cmu;\Q)$ equals
\[
(-1)^{g} \frac{A(t)^{2g}}{n^{4g}\om_{\mu}(t)}
\]
for any $\mu \in \Delta.$ Let $\rho_{r}^{k,\mu}$ denote the
restriction of $\rho_{r}^{k}$ in $\HG(\Cmu;\Q).$ Thus we have that
\[
\frac{(-1)^{g}}{n^{4g}} A(t)^{2g} = \sum_{k=0}^{n-1} \rho_{-1}^{k,\mu}
t^{k} \indent \mbox{mod } \om_{\mu}(t).
\]
Hence by substituting $t = -a_{1}^{i}$ for each $i$ we obtain
\[
\frac{(-1)^{g}}{n^{4g}} \left( \prod_{p=1,p \neq i}^{n} (a_{1}^{i} -
a_{1}^{p}) \right)^{2g} = \sum_{k=0}^{n-1}
\rho_{-1}^{k,\mu}(-a_{1}^{i})^{k}.
\]
Since the direct sum of restriction maps
\[
\HG(\C;\Q) \rightarrow \bigoplus_{\mu \in \Delta} \HG(\Cmu;\Q)
\]
is injective \cite[prop. 3]{K2} we have that
\begin{equation}
\frac{(-1)^{g}}{n^{4g}} \left( \prod_{p=1,p \neq i}^{n} (\delta_{i} -
\delta_{p}) \right)^{2g} = \sum_{k=0}^{n-1}
\rho_{-1}^{k}(-\delta_{i})^{k} \label{321}.
\end{equation}
Solving the equations (\ref{321}) we obtain
\begin{equation}
\rho_{-1}^{k} = \frac{(-1)^{g+n}}{n^{4g}} \sum_{i=1}^{n} S_{i}^{k}
\left( \prod_{p=1,p \neq i}^{n} (\delta_{i} - \delta_{p})
\right)^{2g-1} \label{rev5}
\end{equation}
where $S_{i}^{k}$ equals the $k$th elementary symmetric polynomial in
$\{\delta_{p}: p \neq i\}$.\\
\indent We will show later, in proposition \ref{Pontinad}, that the
relations $\rho_{-1}^{k}$ above are insufficient to prove any
vanishing of the Pontryagin ring in ranks greater than three. For now
consider the rank three case. We write
\[
\alpha = \delta_{1}-\delta_{2},\indent \beta = \delta_{2}-\delta_{3},
\indent \gamma = \delta_{3} - \delta_{1}.
\]
We know from lemma \ref{Pont} that the Pontryagin ring is generated by the
elementary symmetric polynomials in $\alpha^{2},\beta^{2}$ and
$\gamma^{2}.$ The relations $\rho_{-1}^{0},
\rho_{-1}^{1},\rho_{-1}^{2}$ read as
\begin{eqnarray}
(\alpha \beta)^{2g-1} + (\beta \gamma)^{2g-1} + (\gamma \alpha)^{2g-1}
= 0, \label{reva}\\
(\delta_{1}+\delta_{3}) (\alpha \beta)^{2g-1} +(\delta_{2}+\delta_{1})
(\beta \gamma)^{2g-1} + (\delta_{3}+\delta_{2}) (\gamma \alpha)^{2g-1}
= 0, \label{revb}\\
(\delta_{1} \delta_{3})(\alpha \beta)^{2g-1} +(\delta_{2} \delta_{1})
(\beta \gamma)^{2g-1} + (\delta_{3} \delta_{2}) (\gamma \alpha)^{2g-1}
= 0. \label{revc}
\end{eqnarray}
The equations (\ref{reva}), $a_{1} \times$ (\ref{reva}) --
(\ref{revb}), and (\ref{revc}) $+ a_{1} \times$ (\ref{revb}) $-a_{2}
\times$ (\ref{reva}) then show
\begin{equation}
(\delta_{2})^{k} (\alpha \beta)^{2g-1} + (\delta_{3})^{k} (\beta
\gamma)^{2g-1} +(\delta_{1})^{k}  (\gamma \alpha)^{2g-1} = 0, \label{rev2}
\end{equation}
for $k = 0,1,2$. Note that
\[
(\delta_{i})^{r+3} = a_{1} (\delta_{i})^{r+2} - a_{2}
(\delta_{i})^{r+1} + a_{3} (\delta_{i})^{r}
\]
and hence equation (\ref{rev2}) holds for all non-negative
$k$. Further note that
\begin{equation} 
\gamma^{2} = (a_{1})^{2} - 4 a_{2} + 2 a_{1} \delta_{2} - 3
(\delta_{2})^{2} \label{rev1}
\end{equation}
and so combining equation (\ref{rev2}) with equation (\ref{rev1}) and
two similar equations for $\alpha^{2}$ and $\beta^{2}$ we see that
\[
\gamma^{2l} (\delta_{2})^{k} (\alpha \beta)^{2g-1} + \alpha^{2l}
(\delta_{3})^{k} (\beta \gamma)^{2g-1} + \beta^{2l} (\delta_{1})^{k}
(\gamma \alpha)^{2g-1} = 0,
\]
for any non-negative $k,l$. Let $r,s,t$ be three non-negative integers
with an even sum. Note
\[
2 \alpha = (a_{1} - 3 \delta_{2}) - \gamma, \indent
2 \beta =  (3 \delta_{2} - a_{1}) - \gamma,
\]
and hence $(\alpha^{r} \beta^{s} + \alpha^{s} \beta^{r}) \gamma^{t},$ 
when written in terms of $a_{1}, \delta_{2}$ and $\gamma$ is an even
function in $\gamma$.\\
\indent Now any element of the Pontryagin ring can be written
as a sum of elements of the form
\[
F(u,v,w) = \alpha^{u}\beta^{v}\gamma^{w} + \alpha^{v}\beta^{w}\gamma^{u} +
\alpha^{w}\beta^{u}\gamma^{v} + \alpha^{u}\beta^{w}\gamma^{v} +
\alpha^{v}\beta^{u}\gamma^{w} + \alpha^{w}\beta^{v}\gamma^{u},
\]
where $u+v+w$ is even. From the argument above we know that
\begin{equation}
F(2g-1+r,2g-1+s,t) = 0 \label{rev4}
\end{equation}
for $r,s,t \geq 0$ and $r+s+t$ even. If $u \geq 1$ then we have
\begin{equation}
F(u,v,w) = -F(u-1,v,w+1) - F(u-1,v+1,w) \label{rev3}
\end{equation}
since $\alpha+\beta+\gamma=0$.\\
\indent Suppose now that $ u \geq v \geq w$. We
claim $F(u,v,w)=0$ if $u+v+w \geq 6g-4.$ Note that
\[
\mbox{max}\{u,v,w\} > \mbox{max}\{u-1,v+1,w+1\}
\]
unless $u-v$ equals zero or one. In either case we find that $u
\geq v \geq 2g-1$ and hence $F(u,v,w) = 0$ by (\ref{rev4}). Hence by
repeated applications of identity (\ref{rev3}) we see
that $F(u,v,w) = 0$ when $u+v+w \geq 6g-4$ and so we have:\\[\baselineskip]
{\bf THEOREM 2.} {\em The Pontryagin ring of the moduli space ${\cal
M}(3,d)$ vanishes in degrees $12g-8$ and above.}
\begin{rem}
Theorem 2 falls short of Neeman's conjecture
\cite{NE} which states that the Pontryagin ring of $\mnd$
should vanish in degrees above $2gn^{2}-4g(n-1)+2$. When $n=3$ this gives
$10g+2$.
\end{rem}
\begin{rem}
In the rank two case the relations (\ref{rev5}) show
that
\[
((a_{1})^{2} - 4a_{2})^{g}=0
\]
and that the Pontryagin ring of ${\cal M}(2,d)$ vanishes in degrees
greater than or equal to $4g$, duplicating Kirwan's proof of the
Newstead-Ramanan conjecture.
\end{rem}
\indent To conclude we show now that the relations $\rho_{-1}^{k}$ are
inadequate to show any vanishing of the Pontryagin ring when $n \geq 4$. From equation
(\ref{rev5}) we see that the ideal of the Pontryagin ring is contained
in the ideal generated by the formal expressions
\begin{equation}
\left( \prod_{p=1,p \neq i}^{n} (\delta_{i} - \delta_{p})
\right)^{2g-1} \label{rev6}.
\end{equation}
Let $I$ denote the ideal generated by
the relations (\ref{rev6}) and consider this as an ideal of ${\bf C}[\delta_{1},...,\delta_{n}]$. By Hilbert's Nullstellensatz the radical
$\sqrt{I}$ of $I$ consists of those elements of the Pontryagin ring which vanish on
the intersection of the subspaces given by
\begin{equation}
\prod_{p \neq i} (\delta_{i} - \delta_{p}) = 0, \indent i=1,...,n. \label{Q16}
\end{equation}
We shall consider the even and odd cases for $n$ separately.\\
\indent (i) $n$ is even -- write $n=2m.$ The intersection of the
subspaces (\ref{Q16}) consists of $(2m)!/(2^{m}m!)$ distinct
$m$-dimensional subspaces of ${\bf C}^{n}.$ One of these subspaces is
given by the equations
\begin{equation}
\delta_{2k-1} = \delta_{2k}, \indent k=1,...,m. \label{plane}
\end{equation}
We know from lemma \ref{Pont} that the total Pontryagin class $p(T)$ of
$\mnd$ equals
\[
\prod_{1 \leq k < l \leq n} (1 + (\delta_{k}-\delta_{l})^{2})^{2 \g}
\]
and in the subspace (\ref{plane}) $p(T)$ then equals
\[
\prod_{1 \leq k < l \leq m} (1 + (\delta_{2 k -1} - \delta_{2 l
-1})^{2})^{8 \g}.
\]
In particular we see that none of the Pontryagin classes of $\mnd$
vanish on the subspace (\ref{plane}).\\
\indent (ii) $n$ is odd -- write $n=2m+1$. The intersection of the
subspaces (\ref{Q16}) consists of $(2k+1)!/(3 \cdot 2^{k}(k-1)!)$
distinct $k$-dimensional subspaces of ${\bf C}^{n}$. One of these
subspaces is given by the equations
\begin{equation}
\delta_{1} = \delta_{2} = \delta_{3}, \quad \delta_{2k} = \delta_{2k+1}, \quad k=2,...,m. \label{Plane}
\end{equation}
In the subspace (\ref{Plane}) the total Pontryagin class of $\mnd$
equals
\[
\left( \prod_{2 \leq k \leq m} (1+(\delta_{1}-\delta_{2k})^{2})^{12\g}
\right) \left( \prod_{2 \leq k < l \leq m}
(1+(\delta_{2k}-\delta_{2l})^{2})^{8\g} \right).
\]
In particular we see that none of the Pontryagin classes of $\mnd$
vanish on the subspace (\ref{Plane}).\\
\indent Thus we see that none of the Pontryagin classes $p_{r}(T)$ are
nilpotent modulo the formal relations (\ref{rev6}). Hence: 
\begin{prop}
\label{Pontinad}
For $n \geq 4$ the Pontryagin classes $p_{r}(T) \in
H^{4r}(\mnd;\Q)$ are not nilpotent modulo $\rho_{-1}^{k}$ for $0 \leq
k \leq n-1$. In particular these relations are inadequate to
prove any non-trivial vanishing of the Pontryagin ring. 
\end{prop}

\end{document}